\catcode`\@=11

\font\tenmsa=msam10
\font\sevenmsa=msam7
\font\fivemsa=msam5
\font\tenmsb=msbm10
\font\sevenmsb=msbm7
\font\fivemsb=msbm5
\newfam\msafam
\newfam\msbfam
\textfont\msafam=\tenmsa  \scriptfont\msafam=\sevenmsa
  \scriptscriptfont\msafam=\fivemsa
\textfont\msbfam=\tenmsb  \scriptfont\msbfam=\sevenmsb
  \scriptscriptfont\msbfam=\fivemsb

\def\hexnumber@#1{\ifnum#1<10 \number#1\else
 \ifnum#1=10 A\else\ifnum#1=11 B\else\ifnum#1=12 C\else
 \ifnum#1=13 D\else\ifnum#1=14 E\else\ifnum#1=15 F\fi\fi\fi\fi\fi\fi\fi}

\def\msa@{\hexnumber@\msafam}
\def\msb@{\hexnumber@\msbfam}
\mathchardef\boxdot="2\msa@00
\mathchardef\boxplus="2\msa@01
\mathchardef\boxtimes="2\msa@02
\mathchardef\square="0\msa@03
\mathchardef\blacksquare="0\msa@04
\mathchardef\centerdot="2\msa@05
\mathchardef\lozenge="0\msa@06
\mathchardef\blacklozenge="0\msa@07
\mathchardef\circlearrowright="3\msa@08
\mathchardef\circlearrowleft="3\msa@09
\mathchardef\rightleftharpoons="3\msa@0A
\mathchardef\leftrightharpoons="3\msa@0B
\mathchardef\boxminus="2\msa@0C
\mathchardef\Vdash="3\msa@0D
\mathchardef\Vvdash="3\msa@0E
\mathchardef\vDash="3\msa@0F
\mathchardef\twoheadrightarrow="3\msa@10
\mathchardef\twoheadleftarrow="3\msa@11
\mathchardef\leftleftarrows="3\msa@12
\mathchardef\rightrightarrows="3\msa@13
\mathchardef\upuparrows="3\msa@14
\mathchardef\downdownarrows="3\msa@15
\mathchardef\upharpoonright="3\msa@16

\mathchardef\downharpoonright="3\msa@17
\mathchardef\upharpoonleft="3\msa@18
\mathchardef\downharpoonleft="3\msa@19
\mathchardef\rightarrowtail="3\msa@1A
\mathchardef\leftarrowtail="3\msa@1B
\mathchardef\leftrightarrows="3\msa@1C
\mathchardef\rightleftarrows="3\msa@1D
\mathchardef\Lsh="3\msa@1E
\mathchardef\Rsh="3\msa@1F
\mathchardef\rightsquigarrow="3\msa@20
\mathchardef\leftrightsquigarrow="3\msa@21
\mathchardef\looparrowleft="3\msa@22
\mathchardef\looparrowright="3\msa@23
\mathchardef\circeq="3\msa@24
\mathchardef\succsim="3\msa@25
\mathchardef\gtrsim="3\msa@26
\mathchardef\gtrapprox="3\msa@27
\mathchardef\multimap="3\msa@28
\mathchardef\therefore="3\msa@29
\mathchardef\because="3\msa@2A
\mathchardef\doteqdot="3\msa@2B

\mathchardef\triangleq="3\msa@2C
\mathchardef\precsim="3\msa@2D
\mathchardef\lesssim="3\msa@2E
\mathchardef\lessapprox="3\msa@2F
\mathchardef\eqslantless="3\msa@30
\mathchardef\eqslantgtr="3\msa@31
\mathchardef\curlyeqprec="3\msa@32
\mathchardef\curlyeqsucc="3\msa@33
\mathchardef\preccurlyeq="3\msa@34
\mathchardef\leqq="3\msa@35
\mathchardef\leqslant="3\msa@36
\mathchardef\lessgtr="3\msa@37
\mathchardef\backprime="0\msa@38
\mathchardef\risingdotseq="3\msa@3A
\mathchardef\fallingdotseq="3\msa@3B
\mathchardef\succcurlyeq="3\msa@3C
\mathchardef\geqq="3\msa@3D
\mathchardef\geqslant="3\msa@3E
\mathchardef\gtrless="3\msa@3F
\mathchardef\sqsubset="3\msa@40
\mathchardef\sqsupset="3\msa@41
\mathchardef\trianglerighteq="3\msa@44
\mathchardef\trianglelefteq="3\msa@45
\mathchardef\bigstar="0\msa@46
\mathchardef\between="3\msa@47
\mathchardef\blacktriangledown="0\msa@48
\mathchardef\blacktriangleright="3\msa@49
\mathchardef\blacktriangleleft="3\msa@4A
\mathchardef\blacktriangle="0\msa@4E
\mathchardef\triangledown="0\msa@4F
\mathchardef\eqcirc="3\msa@50
\mathchardef\lesseqgtr="3\msa@51
\mathchardef\gtreqless="3\msa@52
\mathchardef\lesseqqgtr="3\msa@53
\mathchardef\gtreqqless="3\msa@54
\mathchardef\Rrightarrow="3\msa@56
\mathchardef\Lleftarrow="3\msa@57
\mathchardef\veebar="2\msa@59
\mathchardef\barwedge="2\msa@5A
\mathchardef\doublebarwedge="2\msa@5B
\mathchardef\angle="0\msa@5C
\mathchardef\measuredangle="0\msa@5D
\mathchardef\sphericalangle="0\msa@5E
\mathchardef\varpropto="3\msa@5F
\mathchardef\smallsmile="3\msa@60
\mathchardef\smallfrown="3\msa@61
\mathchardef\Subset="3\msa@62
\mathchardef\Supset="3\msa@63
\mathchardef\Cup="2\msa@64

\mathchardef\Cap="2\msa@65

\mathchardef\curlywedge="2\msa@66
\mathchardef\curlyvee="2\msa@67
\mathchardef\leftthreetimes="2\msa@68
\mathchardef\rightthreetimes="2\msa@69
\mathchardef\subseteqq="3\msa@6A
\mathchardef\supseteqq="3\msa@6B
\mathchardef\bumpeq="3\msa@6C
\mathchardef\Bumpeq="3\msa@6D
\mathchardef\lll="3\msa@6E

\mathchardef\ggg="3\msa@6F

\mathchardef\circledS="0\msa@73
\mathchardef\pitchfork="3\msa@74
\mathchardef\dotplus="2\msa@75
\mathchardef\backsim="3\msa@76
\mathchardef\backsimeq="3\msa@77
\mathchardef\complement="0\msa@7B
\mathchardef\intercal="2\msa@7C
\mathchardef\circledcirc="2\msa@7D
\mathchardef\circledast="2\msa@7E
\mathchardef\circleddash="2\msa@7F
\def\ulcorner{\delimiter"4\msa@70\msa@70 }
\def\urcorner{\delimiter"5\msa@71\msa@71 }
\def\llcorner{\delimiter"4\msa@78\msa@78 }
\def\lrcorner{\delimiter"5\msa@79\msa@79 }
\def\yen{\mathhexbox\msa@55 }
\def\checkmark{\mathhexbox\msa@58 }
\def\circledR{\mathhexbox\msa@72 }
\def\maltese{\mathhexbox\msa@7A }
\mathchardef\lvertneqq="3\msb@00
\mathchardef\gvertneqq="3\msb@01
\mathchardef\nleq="3\msb@02
\mathchardef\ngeq="3\msb@03
\mathchardef\nless="3\msb@04
\mathchardef\ngtr="3\msb@05
\mathchardef\nprec="3\msb@06
\mathchardef\nsucc="3\msb@07
\mathchardef\lneqq="3\msb@08
\mathchardef\gneqq="3\msb@09
\mathchardef\nleqslant="3\msb@0A
\mathchardef\ngeqslant="3\msb@0B
\mathchardef\lneq="3\msb@0C
\mathchardef\gneq="3\msb@0D
\mathchardef\npreceq="3\msb@0E
\mathchardef\nsucceq="3\msb@0F
\mathchardef\precnsim="3\msb@10
\mathchardef\succnsim="3\msb@11
\mathchardef\lnsim="3\msb@12
\mathchardef\gnsim="3\msb@13
\mathchardef\nleqq="3\msb@14
\mathchardef\ngeqq="3\msb@15
\mathchardef\precneqq="3\msb@16
\mathchardef\succneqq="3\msb@17
\mathchardef\precnapprox="3\msb@18
\mathchardef\succnapprox="3\msb@19
\mathchardef\lnapprox="3\msb@1A
\mathchardef\gnapprox="3\msb@1B
\mathchardef\nsim="3\msb@1C
\mathchardef\napprox="3\msb@1D
\mathchardef\nsubseteqq="3\msb@22
\mathchardef\nsupseteqq="3\msb@23
\mathchardef\subsetneqq="3\msb@24
\mathchardef\supsetneqq="3\msb@25
\mathchardef\subsetneq="3\msb@28
\mathchardef\supsetneq="3\msb@29
\mathchardef\nsubseteq="3\msb@2A
\mathchardef\nsupseteq="3\msb@2B
\mathchardef\nparallel="3\msb@2C
\mathchardef\nmid="3\msb@2D
\mathchardef\nshortmid="3\msb@2E
\mathchardef\nshortparallel="3\msb@2F
\mathchardef\nvdash="3\msb@30
\mathchardef\nVdash="3\msb@31
\mathchardef\nvDash="3\msb@32
\mathchardef\nVDash="3\msb@33
\mathchardef\ntrianglerighteq="3\msb@34
\mathchardef\ntrianglelefteq="3\msb@35
\mathchardef\ntriangleleft="3\msb@36
\mathchardef\ntriangleright="3\msb@37
\mathchardef\nleftarrow="3\msb@38
\mathchardef\nrightarrow="3\msb@39
\mathchardef\nLeftarrow="3\msb@3A
\mathchardef\nRightarrow="3\msb@3B
\mathchardef\nLeftrightarrow="3\msb@3C
\mathchardef\nleftrightarrow="3\msb@3D
\mathchardef\divideontimes="2\msb@3E
\mathchardef\varnothing="0\msb@3F
\mathchardef\nexists="0\msb@40
\mathchardef\mho="0\msb@66
\mathchardef\thorn="0\msb@67
\mathchardef\beth="0\msb@69
\mathchardef\gimel="0\msb@6A
\mathchardef\daleth="0\msb@6B
\mathchardef\lessdot="3\msb@6C
\mathchardef\gtrdot="3\msb@6D
\mathchardef\ltimes="2\msb@6E
\mathchardef\rtimes="2\msb@6F
\mathchardef\shortmid="3\msb@70
\mathchardef\shortparallel="3\msb@71
\mathchardef\smallsetminus="2\msb@72
\mathchardef\thicksim="3\msb@73
\mathchardef\thickapprox="3\msb@74
\mathchardef\approxeq="3\msb@75
\mathchardef\succapprox="3\msb@76
\mathchardef\precapprox="3\msb@77
\mathchardef\curvearrowleft="3\msb@78
\mathchardef\curvearrowright="3\msb@79
\mathchardef\digamma="0\msb@7A
\mathchardef\varkappa="0\msb@7B
\mathchardef\hslash="0\msb@7D
\mathchardef\hbar="0\msb@7E
\mathchardef\backepsilon="3\msb@7F
\def\Bbb{\ifmmode\let\next\Bbb@\else
 \def\next{\errmessage{Use \string\Bbb\space only in math mode}}\fi\next}
\def\Bbb@#1{{\Bbb@@{#1}}}
\def\Bbb@@#1{\fam\msbfam#1}

\catcode`\@=\active

\def\lform{\hbox{$\sqcup$}\llap{\hbox{$\sqcap$}}}
\def\eps{{\epsilon}}

\def\<{\langle}
\def\>{\rangle}

\def\tens{\mathop{\otimes}}
\def\la{{\triangleright}}\def\ra{{\triangleleft}}

\def\proof{\goodbreak\noindent{\bf Proof\quad}}
\def\endproof{{\ $\lform$}\bigskip }

\def\text#1{{\rm #1}}
\def\note#1{}

\documentstyle[11pt,epsfig]{article}
\newtheorem{lemma}{Lemma}[section]
\newtheorem{propos}[lemma]{Proposition}
\newtheorem{example}[lemma]{Example}
\newtheorem{theorem}[lemma]{Theorem}

\newtheorem{defin}[lemma]{Definition}

\begin{document}
\baselineskip 22pt

\newpage
\baselineskip 22pt

{\ }\hskip 4.7in   Swan/00
\vspace{.2in}

\begin{center} {\Large MAKING NON-TRIVIALLY ASSOCIATED TENSOR CATEGORIES 
 FROM LEFT COSET REPRESENTATIVES}
\\ \baselineskip 13pt{\ }
{\ }\\
Edwin Beggs \\
{\ }\\
Department of Mathematics\\
University of Wales, Swansea\\
Singleton Park, Swansea SA2 8PP, UK
\end{center}

\begin{center}
16th February 2000 
\end{center}
\vspace{10pt}
\begin{quote}\baselineskip 13pt
\noindent{\bf ABSTRACT} The paper begins by giving an algebraic structure
on a set of coset representatives for the left action of a subgroup on a group. 
From this we construct a non-trivially associated tensor category. Also a double construction
is given, and this allows the construction of a
non-trivially associated braided tensor category. In this category
we explicitly reconstruct a braided Hopf algebra, whose representations
comprise the category itself. 
\end{quote}
\baselineskip 22pt

\section{INTRODUCTION}
This paper is in some respects a sequel to \cite{BGM:fin}, 
which considered group doublecross products, i.e.\ a group factoring into two
subgroups. 
Group doublecross products are 
the foundation of one way to look at certain integrable
field theories \cite{BJ:form}. Here the space-time
is imbedded in the group by a function called the `classical vacuum map'. This imbedding
possibly encodes information about the geometry of the space-time, and in general there is no reason
why its image should be a subgroup
(also see \cite{Ma:mec}). This raised the possibility of considering more general
factorisations of groups, and their corresponding algebras.
These algebras turned out to be non-trivially
associated,  after the manner of  \cite{AM:oct}. 
In retrospect it is suprising how many of the results of the standard theory of doublecross
products and bicrossproducts
\cite{Tak:mat,Ma:phy} carry over into the present setting. Typically all that is required is the insertion
of a few additional pieces into the relevant formulae. Throught the paper I have liberally used
 \cite{Ma:book} as a reference for tensor categories and braided Hopf algebras. 
The reader should note that the coset representative constructions in this paper are
essentially those appearing in the group cohomology analysis of exact sequences
of groups
\cite{Weiss:book}. This then leads on to non-commutative topological
cohomology via crossed modules \cite{RB:gpcross,{Y:top}}.  It is not clear to me 
whether there is a more direct link between non-trivially associated tensor
categories and cohomology. 

Take a representative element of every coset for the left action of a subgroup $G$ on
a group $X$, and form a set $M$ of these elements. From the algebraic structure on $X$ we can 
construct a binary operation on $M$ which has a left identity and right division. Conversely any
set with such a binary operation can be realised as a set of left coset representatives
for the quotient of two permutation groups. 
This binary operation is
not associative, but the breakdown of associativity is given by a `cocycle' $\tau:M\times M\to G$. 
Using this cocycle we can construct a non-trivial associator for a category ${\cal C}$ of $M$-graded
right representations of $G$. This category also has evaluation and coevaluation maps,
making it into a rigid tensor category. If we assume that the binary operation on $M$
has left division (which is not always true) then the grading and group action
can be combined into the action of an algebra $H$ on the objects in the category. 
It turns out that $H$ itself is in ${\cal C}$, and that the multiplication is
associative (using the associator). 

Next we construct a double factorisation from two copies of the original group $X$, though we give
one copy of $X$ a different binary operation. The category ${\cal D}$ constructed from the
double (in the same way that ${\cal C}$ was constructed from $X$) is braided,
as well as being non-trivially associated. We can again form an algebra $D$ in ${\cal D}$ whose action
combines both the gradings and actions in the definition of ${\cal D}$. But now, using the
braiding in ${\cal D}$, we can reconstruct a coproduct on $D$ from the tensor 
product structure in ${\cal D}$. The existence of such a braided Hopf algebra is guaranteed
by a general reconsruction theory \cite{Ma:book},
and we explicitly calculate the braided Hopf algebra structure on $D$ using these methods. 
I do not know whether the braided category ${\cal D}$ gives
any interesting knot invariants. 

Except for the section on left coset representatives,
I shall assume that all groups are finite and
that all vector spaces are finite dimensional. 
This is to avoid problems with measurability and continuity. 

I would like to thank Y.~Bespalov, Ronnie Brown and Shahn Majid for their assistance
in the preparation of this paper.

%\section*{Preliminaries}

\section{Left coset representatives.}

\begin{defin} Given a group $X$ and a subgroup $G$, call $M\subset X$ a set of left
coset representatives if for every $x\in X$ there is a unique $s\in M$ so that
$x\in Gs$. We shall call the decomposition $x=us$ for $u\in G$ and $s\in M$ the 
unique factorisation of $x$. 
\end{defin}

For the remainder of this section we assume that $M\subset X$ is a set of left
coset representatives for the subgroup $G\subset X$. The identity in $X$
will be denoted $e$.

\begin{defin} Given $s,t\in M$, we define $\tau(s,t)\in G$ and $s\cdot t\in M$
by the unique factorisation $st=\tau(s,t)(s\cdot t)$ in $X$. We also define functions
$\la:M\times G\to G$ and $\ra:M\times G\to M$ by the unique factorisation
$su=(s\la u)(s\ra u)$ for $s,s\ra u\in M$ and $u,s\la u\in G$. 
\end{defin}

\begin{propos} The binary operation $(M,\cdot)$ has a unique left identity $e_m\in M$
(i.e.\ $e_m \cdot t=t$ for all $t\in M$) and has the right division property 
(i.e.\  for all $t,s\in M$ there is a unique solution $p\in M$ to the equation $p\cdot s=t$).
If $e\in M$ then $e_m=e$ is also a right identity.
\end{propos}
\proof There is a unique factorisation $e=ue_m$ for $e_m\in M$ and $u\in G$, 
so $G\cap M=\{e_m\}$. Then $e_m\cdot t=t$ by definition. Conversely if $s\cdot t=t$
then $st=\tau(s,t)t$, so $s\in G\cap M$.

If $p\cdot s=t$ then $ps=\tau(p,s)t$, so $\tau(p,s)^{-1}p=ts^{-1}$. Now apply unique factorisation
to $ts^{-1}\in X$. 
\endproof

By applying the right division property to solve the equation $p\cdot t=e_m$ for a given $t\in M$, 
we see that there is a unique left inverse $t^L$ for every $t$, satisfying the equation
$t^L\cdot t=e_m$. 

We shall use the result of the next proposition at many places in the paper:

\begin{propos}\label{rels} The following identities between $(M,\cdot)$ and $\tau$ hold,
where we take $t,s,p\in M$ and $u,v\in G$:
\begin{eqnarray*}
s\la(t\la u)\ =\ \tau(s,t)\, ((s\cdot t)\la u)\, \tau(s\ra (t\la u),t\ra u)^{-1}
 &{\rm and}& (s\cdot t)\ra u \ =\ (s\ra (t\la u) )\cdot (t\ra u)
\ ,  \cr
s\la uv\ =\ (s\la u)\, ((s\ra u)\la v)  &{\rm and}&
s\ra uv\ =\ (s\ra u)\ra v \ ,\cr
\tau(p,s)\tau(p\cdot s,t)\ =\ (p\la \tau(s,t))\, \tau(p\ra \tau(s,t),s\cdot t)
&{\rm and}&
(p\ra \tau(s,t))\cdot(s\cdot t)\ =\ (p\cdot s)\cdot t\ .
\end{eqnarray*}
\end{propos} 
\proof
We can deduce these identities from the associativity
of $X$. From the identity $(st)u=s(tu)$ we can deduce that
\begin{eqnarray*}
(st)u &=& \tau(s,t)(s\cdot t)u\ =\ \tau(s,t)((s\cdot t)\la u)((s\cdot t)\ra u)\ ,\cr
s(tu) &=& s(t\la u)(t\ra u)\ =\ (s\la(t\la u))(s\ra(t\la u))(t\ra u)   \cr &=&
(s\la(t\la u))\, \tau(s\ra(t\la u),t\ra u)\, ((s\ra(t\la u))\cdot(t\ra u))\ .
\end{eqnarray*}
The first line follows by uniqueness of factorisation.
From $s(uv)=(su)v$,
\begin{eqnarray*}
s(uv) &=& (s\la uv)\, (s\ra uv)\ ,\cr
(su)v &=& (s\la u)\, (s\ra u)v\ =\ (s\la u)\, ((s\ra u)\la v)\,  ((s\ra u)\ra v)\ ,
\end{eqnarray*}
giving the second line identities. Finally from $p(st)=(ps)t$,
\begin{eqnarray*}
p(st) &=& p\, \tau(s,t)\, (s\cdot t)\ =\ (p\la \tau(s,t))\,  (p\ra \tau(s,t))\, (s\cdot t) \cr
&=& (p\la \tau(s,t))\,  \tau(p\ra \tau(s,t),s\cdot t)\, ((p\ra \tau(s,t))\cdot (s\cdot t))\ ,\cr
(ps)t &=& \tau(p,s)\, (p\cdot s)\, t\ =\ \tau(p,s)\, \tau(p\cdot s,t)\ 
((p\cdot s)\cdot t)\ ,
\end{eqnarray*}
giving the last line.
\endproof

\begin{propos}\label{rels2} The following identities between $(M,\cdot)$ and $\tau$ hold,
for all $t\in M$ and $v\in G$:
\begin{eqnarray*}
e_m\ra v \ =\ e_m\ ,\quad e_m\la v \ =\  e_m\,v\,e_m^{-1}\ &,&\quad t\la e\ =\ e\ ,\quad
t\ra e\ =\ t\ ,\cr \tau(e_m,t)\ =\  e_m\ ,\quad
t\la e_m^{-1}&=& \tau(t\ra e_m^{-1},e_m)^{-1}\ ,\quad 
(t\ra e_m^{-1})\cdot e_m\ =\ t\ .
\end{eqnarray*}
\end{propos} 
\proof We have the factorisation $e_m\, v=(e_m\,v\,e_m^{-1})e_m$, where
$e_m\in M$ and $e_m\,v\,e_m^{-1}\in G$. Also $te=et$ for $e\in G$ and $t\in M$. Next
$e_m t=\tau(e_m,t)(e_m\cdot t)=\tau(e_m,t)t$. Finally
\[
t\ =\ te_m^{-1}e_m\ =\ (t\la e_m^{-1})(t\ra e_m^{-1})e_m\ =\ (t\la e_m^{-1})\, 
\tau(t\ra e_m^{-1}, e_m)\, ((t\ra e_m^{-1})\cdot e_m)\ ,
\]
giving the last identities.\endproof 

This last proposition makes sense because $e_m\in G\cap M$. For situations where
it is convenient to forget about the original group $X$, and just concentrate on $G$ and $(M,\cdot)$,
we will use $e_m\in M$ for the left identity in $(M,\cdot)$, and set $f_m=e_m\in G$.

\begin{example}\label{ex1} Take $X$ to be the permutation group $S_3$ of 3 objects $\{1,2,3\}$,
and let $G$ be the non-normal subgroup $\{e,(12)\}$. Take the
 set of left
coset representatives 
$M=\{(12),(13),(23)\}$.
The dot and  $\tau$ operation are given by the following tables, where the 
row $s$ column $t$ entry corresponds to $s\cdot t$ or $\tau(s,t)$:\newline
 \begin{tabular}{c|ccc}
 $\cdot$  & (12) & (13) & (23) \\
\hline 
(12)   & (12) & (13) & (23) \\
(13)   & (23) & (12) & (13) \\
(23)   & (13) & (23) & (12) \\
\end{tabular}  \qquad 
 \begin{tabular}{c|ccc}
$\tau$   & (12) & (13) & (23) \\
\hline 
(12)   & (12) &(12) & (12) \\
(13)   & (12) &(12) & (12) \\
(23)   & (12) & (12) & (12) \\
\end{tabular} \newline
The fact that $(M,\cdot)$ satisfies the right division property
 is just the condition that every element of $M$ appears exactly
once in each column of the table for dot. 
In this case we also see that every element of $M$ appears exactly
once in each row of the $(M,\cdot)$ table, so $(M,\cdot)$
satisfies left division. However there is no
right identity, so $(M,\cdot)$ does not form a group. 
\end{example}

\begin{example}\label{ex2} Take $X$ to be the permutation group $S_3$ of 3 objects $\{1,2,3\}$,
and let $G$ be the non-normal subgroup $\{e,(12)\}$. Take the
 set of left
coset representatives 
$M=\{e,(23),(13)\}$.  The operation
$\la$ is trivial, and $\ra$ is given by the action of $(12)$ on $M$ swapping
$(23)$ and $(13)$. 
The dot and  $\tau$ operation are given by the following tables:
\newline
 \begin{tabular}{c|ccc}
 $\cdot$  & e & (13) & (23) \\
\hline 
e   & e & (13) & (23) \\
(13)   & (13) & e & (13) \\
(23)   & (23) & (23) & e \\
\end{tabular}  \qquad 
 \begin{tabular}{c|ccc}
 $\tau$  & e & (13) & (23) \\
\hline 
e   & e &e & e \\
(13)   & e & e & (12) \\
(23)   & e & (12) & e \\
\end{tabular} \newline
This time we see that $(M,\cdot)$ does not satisfy left division, but does have
a 2-sided identity. 
\end{example}

\begin{example}\label{ex3} Take $X$ to be the permutation group $S_3$ of 3 objects $\{1,2,3\}$,
and let $G$ be the non-normal subgroup $\{e,(12)\}$. Take the
 set of left
coset representatives 
$M=\{e,(123),(132)\}$.  In this case $M$ is a subgroup, and the subgroup operation is the
dot product.  This is just the case of a group doublecross product \cite{BGM:fin}.
\end{example}

\begin{example} Take $X$ to be the dihedral group $D_6=\<x,y:x^6=y^2=e,xy=yx^5\>$,
and $G$ to be the non-abelian normal subgroup of order 6 generated
by $x^2$ and $y$. We choose $M=\{e,x\}$. 
The dot operation on $M$ is given by $e$ the 2-sided identity
and $x\cdot x=e$. 
The $\tau$ function is given by
$\tau(x,x)=x^2$, and all other combinations giving $e$.
The operation $\ra$ is trivial, and $\la$ is given by $x$ acting
on $G=\{e,x^2,x^4,y,yx^2,yx^4\}$ as the permutation
$(y,yx^4,yx^2)$, i.e.\ $x\la y=yx^4$ etc.
Observe that though $(M,\cdot)$ is a group, $\la$ is not a group action. 
\end{example}

\begin{example} Take $X$ to be the group $S_5$ of permutations of the objects $\{1,2,3,4,5\}$, 
and $G$ to be the subgroup fixing the object $1$. We choose
$M=\{e,(12)(354),(14253),(15234),(13245)\}$. If we set $a=(12)(354)$, $b=(14253)$,
$c=(15234)$ and $d=(13245)$, we get the tables:
\newline
 \begin{tabular}{c|ccccc}
 $\cdot$  & e & a & b & c & d \\
\hline 
e   & e & a & b & c & d \\
a   & a & e & c & d & b \\
b   & b & c & d & a & e \\
c   & c & d & e & b & a \\
d   & d & b & a & e & c 
\end{tabular}  \qquad 
 \begin{tabular}{c|ccccc}
 $\tau$  & e & a & b & c & d \\
\hline 
e   & e & e & e & e & e \\
a   & e & (345) & (2534) & (2345) & (2453) \\
b   & e & (34) & (354) & (2345) & (354) \\
c   & e & (45) & (354) & (254) & (2453) \\
d   & e & (35) & (2534) & (354) & (235) 
\end{tabular} \newline
We see that $(M,\cdot)$ satisfies right and left division, that $e$ is a 2-sided
identity, but that $(M,\cdot)$ is not a group.
\end{example}

The last example is just an application of a general construction:

\begin{theorem}: {\bf A modified Cayley's theorem.}\ Any set 
$F$ with a binary operation $\bullet$
which has a left identity and right division can be imbedded in $S_F$
(the group of permutations of the elements of $F$), as a
set of left coset representatives for the subgroup $G\subset S_F$ which  fixes the left
identity.
\end{theorem} 
\proof The function $\sigma:F\to S_F$ is defined by
$\sigma(g)(f\bullet g)=f$. Note that $\sigma(g)$ is a 1-1 correspondence because
$(F,\bullet)$ has right division. 

Let $e_F$ be the left identity in $(F,\bullet)$, which is unique by right division.
Take any $\phi\in  S_F$, and set $g=\phi^{-1}(e_F)$. Then
$\psi=\phi\circ\sigma(g)^{-1}$ has the property that $\psi(e_F)=e_F$, i.e.\ 
$\psi\in G$. Further if $\chi \sigma(g)=\sigma(f)$ for any $\chi\in G$, then
by applying $\sigma(g)^{-1}\chi^{-1}=\sigma(f)^{-1}$ to $e_F$ we see that
$g=f$. We conclude that the image of $\sigma$ forms a set of 
left coset representatives for the subgroup $G$. 

Now consider the equation $\sigma(f)\sigma(g)=\chi \sigma(h)$, for $\chi\in G$. 
Applying the inverse of each side to $e_F$ we see $ \sigma(g)^{-1}(\sigma(f)^{-1}(e_F))=
 \sigma(g)^{-1}(f)=f\bullet g$, so $h=f\bullet g$ as required. 
\endproof

\begin{propos} The subgroup $G$ is normal in $X$ if and only if $\ra$ is trivial, 
i.e.\ $s\ra u=s$ for all $s\in M$ and $u\in G$.
If $G$ is normal then $(M,\cdot)$ is isomorphic to the quotient group $G\backslash X$,
the isomorphism being
 the restriction of the quotient map $X\to G\backslash X$. The subset $M$ is a subgroup
of $X$
 if and only if $e\in M$ and $\tau(s,t)=e$ for all $s,t\in M$. 
\end{propos}
\proof  Since $G$ is closed under conjugation by elements of $G$, we just have to check
conjugation by elements of $M$ to see if $G$ is normal. Then for all $s\in M$ and $u\in G$:
\[
sus^{-1}\ =\ (s\la u)(s\ra u)s^{-1}\in G \Leftrightarrow s\ra u\in Gs \Leftrightarrow s\ra u=s\ .
\]
If $G$ is normal then we just use the usual definition of multiplication of left cosets. 

If $\tau(s,t)=e$ for all $s,t\in M$, then $st=\tau(s,t)(s\cdot t)=s\cdot t\in M$, so
the subset $M$ is closed under multiplication in $X$. If in addition $e\in M$ then for every $t\in M$
(by right division) there is a $t^L\in M$ so that $t^L\cdot t=e$. Then $t^Lt=\tau(t^L,t)(t^L\cdot t)=e$, 
so $t^L=t^{-1}$, and $M$ is closed under inverse in $X$. 
 \endproof

We would like to remove the dependence on the group $X$, and say that certain
conditions on $G$, $(M,\cdot)$, $\tau$, $\ra$
 and $\la$ are equivalent to the
existence of the group $X$. To this end, for the remainder of this section we forget how 
$(M,\cdot)$ was constructed, and just begin with a group $G$ and a set with binary operation $(M,\cdot)$.

\begin{propos} Suppose that the functions $\la:M\times G\to G$, $\ra:M\times G\to M$
and $\tau:M\times M\to G$ satisfy the identities in (\ref{rels}). Then
the binary operation on the set $G\times M$ defined by 
$
(u,s)(v,t)\, =\, (u(s\la v)\,\tau(s\ra v,t),(s\ra v)\cdot t)
$
is associative.
\end{propos} 
\proof Begin by calculating
\begin{eqnarray}
((u,s)(v,t))(w,p) &=& (\, u(s\la v)\,\tau(s\ra v,t)\, ((s\ra v)\cdot t\la w)\, 
\tau(((s\ra v)\cdot t)\ra w,p) \, , \, (((s\ra v)\cdot t)\ra w)\cdot p) \ ,\cr
(u,s)((v,t)(w,p)) &=& (u,s)\, (v(t\la w)\, \tau(t\ra w,p),(t\ra w)\cdot p) \cr
&=& (\, u\, (s\la v(t\la w) \tau(t\ra w,p))\, \tau(s\ra v(t\la w) \tau(t\ra w,p),
(t\ra w)\cdot p) \cr &&\qquad\, ,\, (s\ra v(t\la w) \tau(t\ra w,p))\cdot ((t\ra w)\cdot p)\, )\ .
\label{ass1}
\end{eqnarray}
To show that the $M$ components of (\ref{ass1}) are identical, we use
\[
(((s\ra v)\cdot t)\ra w)\cdot p \ =\ ((s\ra v(t\la w))\cdot (t\ra w))\cdot p
\ =\ (s\ra v(t\la w) \tau(t\ra w,p))\cdot ((t\ra w)\cdot p)\ .
\]
We can use the identities to show that the $G$ component of $((u,s)(v,t))(w,p)$
is
\[
u(s\la v)\, ((s\ra v)\la(t\la w))\, \tau(s\ra v(t\la w),t\ra w)\, 
\tau((s\ra v(t\la w))\cdot(t\ra w),p)\ ,
\]
wheras the $G$ component of $(u,s)((v,t)(w,p))$ is
\[
u(s\la v)\, ((s\ra v)\la(t\la w))\, ((s\ra v(t\la w))\la \tau(t\ra w,p))\, 
\tau(s\ra v(t\la w) \tau(t\ra w,p),
(t\ra w)\cdot p)\ ,
\]
and then use the identities again to show that these are the same.\endproof

\begin{propos}\label{inv} Suppose that the functions $\la:M\times G\to G$, $\ra:M\times G\to M$
and $\tau:M\times M\to G$ satisfy the identities in (\ref{rels}). Suppose that there is a left identity
$e_m\in M$ for $(M,\cdot)$ and an element $f_m\in G$
so that for all $t\in M$ and $v\in G$,
\begin{eqnarray*}
e_m\ra v \ =\ e_m\ ,\quad e_m\la v \ =\  f_m\,v\,f_m^{-1}\ &,&\quad t\la e\ =\ e\ ,\quad
t\ra e\ =\ t\ ,\cr \tau(e_m,t)\ =\  f_m\ ,\quad
t\la f_m^{-1}&=& \tau(t\ra f_m^{-1},e_m)^{-1}\ ,\quad 
(t\ra f_m^{-1})\cdot e_m\ =\ t\ .
\end{eqnarray*}
Then the multiplication on the set $G\times M$ defined in the previous proposition
has a 2-sided identity $(f_m^{-1},e_m)$. 

If in addition we suppose that $(M,\cdot)$ has left inverses (i.e.\ for every $t\in M$ there
is a $t^L\in M$ so that $t^L\cdot t=e_m$), then $G\times M$ has left inverses, defined by
\[
(v,t)^{L}\ =\ (f_m^{-1}\tau(t^L,t)^{-1}(t^L\la v^{-1}),t^L\ra v^{-1})\ .
\]
These properties imply that $G\times M$ with the given structure is a group.
\end{propos} 
\proof First we check the 2-sided identity
\begin{eqnarray*}
(f_m^{-1},e_m)(v,t) &=& (f_m^{-1}(e_m\la v)\, \tau(e_m\ra v,t),(e_m\ra v)\cdot t) \cr
&=& (f_m^{-1}f_m v f_m^{-1}f_m, e_m\cdot t)\ =\ (v,t)\ ,  \cr
(v,t)(f_m^{-1},e_m) &=& (v(t\la f_m^{-1})\tau(t\ra f_m^{-1},e_m),(t\ra f_m^{-1})\cdot e_m)
\ =\ (v,t)\ .
\end{eqnarray*}
Finally we check the left inverse:
\begin{eqnarray*}
(f_m^{-1}\tau(t^L,t)^{-1}(t^L\la v^{-1}),t^L\ra v^{-1})(v,t)&=&
(\, f_m^{-1}\tau(t^L,t)^{-1}(t^L\la v^{-1})\, ((t^L\ra v^{-1})\la v)\, \tau(t^L,t)\, ,\, t^L\cdot t\, )\cr
&=& (\, f_m^{-1}\tau(t^L,t)^{-1}(t^L\la v^{-1}v)\, \tau(t^L,t)\, ,\,e_m\, )\ =\ (f_m^{-1},e_m)\ .
\end{eqnarray*}
It is now standard algebra to check that these conditions on identities and inverses,
together with associativity, give a group structure.\endproof

We can now imbed the group $G$ in $G\times M$ by the map $v\mapsto (vf_m^{-1},e_m)$, and
$M$ in $G\times M$ by the map $t\mapsto (e,t)$. Then we get the original situation with
left coset representatives.

\section{A tensor category}
Take a group $X$ with subgroup $G$, and
a set of left coset representatives $M$. We again take $e_m$ to be the left identity in $M$
and $f_m$ to be the corresponding element in $G$. If the reader wishes, the situation
can be simplified if $e\in M$, as then $e_m=f_m=e$. 

Take a category $\cal C$ of finite dimensional vector spaces
over a field $k$, whose objects are right representations of the group $G$ and possess $M$-gradings, 
i.e.\ an object $V$ decomposes as a direct sum of subspaces
 $V=\oplus_{s\in M} V_s$.  If $\xi\in V_s$ for some $s\in M$ we say that
$\xi$ is a {\it homogenous} element of $V$, with grade 
$\langle \xi\rangle=s$. In our formulae in this paper we shall usually assume that
we have chosen homogenous elements of the relevant objects, as the general elements
are just linear combinations of the homogenous elements. 
 We write the action for the representation as $\bar\ra:V\times G\to V$. 
In addition we suppose that the action and grading satisfy a compatibility condition, 
$\langle \xi\bar\ra u\rangle=\langle \xi\rangle\ra u$. 
 The morphisms are linear maps which preserve both the grading and the action, 
i.e.\ for a morphism $\theta:V\to \tilde V$ we have $\langle \theta(\xi)\rangle =\langle \xi\rangle $
and $\theta(\xi)\bar\ra u=\theta(\xi\bar\ra u)$ for all $\xi\in V$ and $u\in G$. 

\begin{propos} 
We can make ${\cal C}$ into a tensor category by taking $V\tens W$ to be the usual vector space
tensor product, with actions and gradings given by
\[
\langle \xi\tens \eta\rangle =\langle \xi\rangle \cdot\langle  \eta\rangle \quad{\rm and}\quad
(\xi\otimes\eta)\bar\ra u=\xi\bar\ra(\langle \eta\rangle \la u)\otimes \eta\bar\ra u\ .
\]
For morphisms $\theta:V\to \tilde V$ and $\phi:W\to \tilde W$ we define
the morphism $\theta\tens\phi:V\tens W\to \tilde V\tens \tilde W$ by $(\theta\tens\phi)(\xi\tens\eta)=
\theta(\xi)\tens\phi(\eta)$, which is just the usual vector space formula.
\end{propos}
\proof  We must check that $\langle (\xi\tens \eta)\bar\ra u\rangle =\langle \xi\tens \eta\rangle \ra u$, which
is automatic from the usual identities (\ref{rels}). Also we have to check that
$((\xi\tens \eta)\bar\ra u)\bar \ra v=(\xi\tens \eta)\bar\ra uv$, which is again simple 
from the identities. It is also straightforward to check that $\theta\tens\phi$ is a morphism in 
${\cal C}$. 
 \endproof

\begin{propos} 
The identity for the tensor operation is just the
vector space $k$ with trivial $G$-action and grade $e_m\in M$. For any object $V$ the morphisms
$l_V:V\to V\tens k$ and $r_V:V\to k\tens V$ are given by the formulae
$l_V(\xi)=\xi\bar\ra f_m^{-1}\tens 1$ and $r_V(\xi)=1\tens \xi$, where $1$ is the 
multiplicative identity in $k$. 
\end{propos}
\proof  We must check that the maps $l_V$ and $r_V$ are morphisms in ${\cal C}$. We have
\begin{eqnarray*}
(\xi\bar\ra f_m^{-1}\tens 1)\bar\ra u &=& \xi\bar\ra f_m^{-1}(\langle 1\rangle \la u)\tens 1\bar\ra u\ =\ 
\xi\bar\ra f_m^{-1}(e_m\la u)\tens 1\ =\ 
\xi\bar\ra uf_m^{-1}\tens 1\ ,\cr
(1\tens \xi)\bar\ra u &=& 1\bar\ra (\langle \xi\rangle \la u)\tens \xi\bar\ra u\ =\ 1\tens
\xi\bar\ra u\ .
\end{eqnarray*}
For the grades we note that $(\langle \xi\rangle \bar\ra f_m^{-1})\cdot e_m=
\langle \xi\rangle $ and $e_m\cdot \langle \xi\rangle =\langle \xi\rangle $, using (\ref{rels2}). 
\endproof

\begin{propos} There is an associator $\Phi_{UVW}:(U\otimes V)\otimes W\to U\otimes (V\otimes W)$
given by
$$
\Phi((\xi\otimes\eta)\otimes\zeta)\ =\ \xi\bar\ra\tau(\langle \eta\rangle ,\langle \zeta\rangle )
\otimes(\eta\otimes\zeta)\ .
$$
\end{propos}
\proof First we must check that $\Phi$ preserves the grading. This is just the identity
\[
(\langle \xi\rangle \cdot\langle \eta\rangle )\cdot\langle \zeta\rangle \ =\ 
(\langle \xi\rangle \ra\tau(\langle \eta\rangle ,\langle \zeta\rangle ))
\cdot(\langle \eta\rangle \cdot \langle \zeta\rangle )\ .
\]
Now we check that the $G$ action
 commutes with the associator. Begin with $\bar\ra u$
$$
\Big(\Big(\xi\otimes\eta\Big)\otimes\zeta\Big)\bar\ra u\ =\ 
\Big(\xi\bar\ra(\langle \eta\rangle \la(\langle \zeta\rangle \la u))
\otimes\eta\bar\ra(\langle \zeta\rangle \la u)\Big)\otimes\zeta\bar\ra u\ ,
$$
and apply $\Phi$ to get
$$
\xi\bar\ra(\langle \eta\rangle \la(\langle \zeta\rangle \la u))\tau(
\langle \eta\rangle \ra(\langle \zeta\rangle \la u), \langle \zeta\rangle \ra u)
\otimes\Big(\eta\bar\ra(\langle \zeta\rangle \la u)\otimes\zeta\bar\ra u\Big)\ .
$$
Applying $\Phi$ first and then $\bar\ra u$ we get
$$
\Big(\xi\bar\ra\tau(\langle \eta\rangle ,\langle \zeta\rangle )
\otimes\Big(\eta\otimes\zeta\Big)\Big)\bar\ra u\ =\ 
\xi\bar\ra\tau(\langle \eta\rangle ,\langle \zeta\rangle )
(\langle \eta\rangle \cdot\langle \zeta\rangle \la u)\otimes
\Big(\eta\bar\ra(\langle \zeta\rangle \la u)\otimes\zeta\bar\ra u\Big)\ ,
$$
which is identical to the first expression by
the usual identities.

Now we must check that $\Phi$ obeys the pentagon condition, which states that
the following two re-bracketings are the same:
\begin{eqnarray*}
((V\tens W)\tens Z)\tens U 
&\to& (V\tens W)\tens (Z\tens U)
\to V\tens (W\tens (Z\tens U))   \cr
((V\tens W)\tens Z)\tens U 
&\to& (V\tens (W\tens Z))\tens U 
\to V\tens ((W\tens Z)\tens U )
\to V\tens (W\tens (Z\tens U)) 
\end{eqnarray*}
We apply these operations to $((\xi\tens\eta)\tens\zeta)\tens\upsilon$, giving
\begin{eqnarray*}
((\xi\tens\eta)\tens\zeta)\tens\upsilon &\mapsto&
(\xi\tens\eta)\bar\ra\tau(\langle \zeta\rangle ,\langle \upsilon\rangle )\tens(\zeta\tens\upsilon) 
\cr &=& 
(\xi\bar\ra(\langle \eta\rangle \la \tau(\langle \zeta\rangle ,\langle \upsilon\rangle ))
\tens\eta\bar\ra\tau(\langle \zeta\rangle ,\langle \upsilon\rangle ))\tens(\zeta\tens\upsilon)
\cr &\mapsto&
\xi\bar\ra(\langle \eta\rangle \la \tau(\langle \zeta\rangle ,\langle \upsilon\rangle ))
\tau(\langle \eta\rangle \ra\tau(\langle \zeta\rangle ,\langle \upsilon\rangle ),
\langle \zeta\rangle \cdot\langle \upsilon\rangle )  \cr &&\quad
\tens(\eta\bar\ra\tau(\langle \zeta\rangle ,\langle \upsilon\rangle )\tens(\zeta\tens\upsilon))\ ,\cr
((\xi\tens\eta)\tens\zeta)\tens\upsilon &\mapsto&
(\xi\bar\ra\tau(\langle \eta\rangle ,\langle \zeta\rangle )\tens(\eta\tens\zeta))\tens\upsilon
\cr  &\mapsto&
\xi\bar\ra\tau(\langle \eta\rangle ,\langle \zeta\rangle )
\tau(\langle \eta\rangle \cdot \langle \zeta\rangle ,\langle \upsilon\rangle )
\tens((\eta\tens\zeta)\tens\upsilon)
\cr  &\mapsto&
\xi\bar\ra\tau(\langle \eta\rangle ,\langle \zeta\rangle )
\tau(\langle \eta\rangle \cdot \langle \zeta\rangle ,\langle \upsilon\rangle )
\tens(\eta\bar\ra\tau(\langle \zeta\rangle ,\langle \upsilon\rangle )\tens(\zeta\tens\upsilon))\ .
\end{eqnarray*}
These are the same by the usual identities. 

We must check the triangle identity, that is the maps
${\rm id}\tens r$ and $\Phi\circ(l\tens {\rm id}):V\tens W\to V\tens(k\tens W)$ 
are the same. 
\begin{eqnarray*}
({\rm id}\tens r)(\xi\tens\eta) &=& \xi\tens(1\tens\eta) \cr 
\Phi\circ(l\tens {\rm id})(\xi\tens\eta) &=& \Phi((\xi\bar\ra f_m^{-1}\tens 1)\tens\eta)\ =\ 
\xi\bar\ra f_m^{-1}\tau(\langle 1\rangle ,\langle \eta\rangle )\tens (1\tens\eta)\ .
\end{eqnarray*}
These are the same as $\tau(\langle 1\rangle ,\langle \eta\rangle )=\tau(e_m,\langle \eta\rangle )=f_m$,
from (\ref{rels2}).

Finally we check that condition that $\Phi$ is a natural transformation, i.e.\ that
the following diagram commutes,
\[\begin{array}{ccc}
(U\tens V)\tens W & \stackrel{\Phi_{UVW}}\longrightarrow & U\tens (V\tens W) \\
\downarrow\ (\psi\tens\theta)\tens\phi & & \downarrow\ \psi\tens(\theta\tens\phi) \\
(\tilde U\tens \tilde V)\tens \tilde W & \stackrel{
\Phi_{\tilde U\tilde V\tilde W}}\longrightarrow & \tilde U\tens (\tilde V\tens \tilde W) 
\end{array} \ ,\]
where
\begin{eqnarray*}
((\psi\tens\theta)\tens\phi)((\xi\tens\eta)\tens\kappa)&=&(\psi(\xi)\tens\theta(\eta))\tens
\phi(\kappa)\ ,\cr
(\psi\tens(\theta\tens\phi))(\xi\tens(\eta\tens\kappa))&=&
\psi(\xi)\tens(\theta(\eta)\tens
\phi(\kappa))\ .
\end{eqnarray*}
 This is simple to check, remembering that the morphisms preserve the grade
and action. 
\endproof

\section{A rigid tensor category}
Take a group $X$ with subgroup $G$, and
a set of left coset representatives $M$ which contains $e$. 
We suppose that $(M,\cdot)$ has right inverses, i.e.\ for every
$s\in M$ there is an $s^R\in M$ so that $s\cdot s^R=e$. 

Take a decomposition of an object $V$ in ${\cal C}$ according to the
grading, i.e.\ $V=\oplus_{s\in M} V_s$, where $\xi\in V_s$ corresponds to 
$\langle \xi\rangle=s$. Now take the dual
vector space $V'$, and set
\[
V'_{s^L}\ =\ \{\alpha\in V' : \alpha|_{V_t}=0\quad\forall t\neq s\}\ .
\]
Then $V'=\oplus_{s\in M} V'_{s^L}$, and we define $\langle \alpha\rangle=s^L$ when
$\alpha\in V'_{s^L}$. The evaluation map ${\rm ev}:V'\tens V\to k$ is defined by
${\rm ev}(\alpha,\xi)=\alpha(\xi)$. We have designed the grading on $V'$ so that
this map preserves gradings. Now considering the action $\bar\ra u$, 
 if we apply evaluation to  $\alpha\bar\ra(\langle \xi\rangle\la u)\tens \xi\bar\ra u$ we should
 get $\alpha(\xi)\bar\ra u=\alpha(\xi)$. To do this we define
$
(\alpha\bar\ra(\langle \xi\rangle\la u))\, ( \xi\bar\ra u)= \alpha(\xi)
$,
or if we put $\eta=\xi\bar\ra u$, 
\[
(\alpha\bar\ra((\langle \eta\rangle\ra u^{-1})\la u))\, ( \eta)\ =\ \alpha(\eta\bar\ra u^{-1})\ =\ 
(\alpha\bar\ra(\langle \eta\rangle\la u^{-1})^{-1})\, ( \eta)\ .
\]
If we rearrange this to give $\alpha\ra v$ we get the following formula;
\begin{eqnarray}
(\alpha\bar\ra v)(\eta)\ =\ \alpha(\eta\bar\ra
\tau(\langle \eta\rangle^L,\langle \eta\rangle)^{-1}\ (\langle \eta\rangle^L\la v^{-1})\ 
\tau(\langle \eta\rangle^L\ra v^{-1},(\langle \eta\rangle^L\ra v^{-1})^R))\  .\label{dact}
\end{eqnarray}

To define the coevaluation map we take a basis $\{\xi\}$ of each $V_s$, and a corresponding
dual basis $\{\hat\xi\}$ of each $V'_{s^L}$, i.e.\ $\hat\eta(\xi)=\delta_{\xi,\eta}$. 
Then we put these bases together for all $s\in M$, and define 
\[
{\rm coev}(1)\ =\ \sum_{\xi\in{\rm basis}}
 \xi\bar\ra \tau(\langle \xi\rangle^L,\langle \xi\rangle)^{-1}\tens\hat\xi\ .
\]

\begin{propos} The coevaluation map defined above is a morphism in ${\cal C}$.
\end{propos}
\proof
First show that each summand in the coevaluation has grade $e$. 
If we put $s=\langle \xi\rangle$, we have to show that $(s\ra\tau(s^L,s)^{-1})\cdot s^L=e$. If we apply
$\cdot s$ to $(s\ra\tau(s^L,s))\cdot s^L$ we get
\[
((s\ra\tau(s^L,s)^{-1})\cdot s^L)\cdot s\ =\ s\cdot (s^L\cdot s)\ =\ s\cdot e\ =\ s\ ,
\]
so using right division shows that $(s\ra\tau(s^L,s)^{-1})\cdot s^L=e$ as required.

 It is reasonably easy to see that the map is independent of the choice of basis. If we apply $\bar\ra u$
to the coevaluation, we get
\[
{\rm coev}(1)\bar\ra u\ =\ 
\sum_\xi \xi\bar\ra \tau(\langle \xi\rangle^L,\langle \xi\rangle)^{-1}(
\langle \xi\rangle^L\la u)
\tens\hat\xi\bar\ra u\ .
\]
Now define a new basis by $\eta=\xi\bar\ra\tau(\langle \xi\rangle^L,\langle \xi\rangle)^{-1}(
\langle \xi\rangle^L\la u)\tau(\langle \xi\rangle^L\ra u,(\langle \xi\rangle^L\ra u)^R)$.
We see that $(\alpha\bar\ra u^{-1})(\xi)=\alpha(\eta)$, so the dual basis is given by
$\hat\xi=\hat\eta\bar\ra u^{-1}$. Now if we write the coevaluation in terms of the new basis
we get
\[
{\rm coev}(1)\ =\ \sum_\eta \eta\bar\ra \tau(\langle \eta\rangle^L,\langle \eta\rangle)^{-1}\tens\hat\eta\ .
\]
Since $\langle\hat \eta\rangle=\langle \hat\xi\rangle\ra u$, we see
that $\tau(\langle \eta\rangle^L,\langle\eta\rangle)=
\tau(\langle \xi\rangle^L\ra u,(\langle \xi\rangle^L\ra u)^R)$, so the expressions for 
${\rm coev}(1)\bar\ra u$ and ${\rm coev}(1)$ in the new basis coincide. We conclude that the action is trivial
on ${\rm coev}(1)$ as required.\endproof

Now we need to check the consistency of the evaluation, coevaluation and associator.
Consider the maps, for a homogenous basis element $\eta$:
\begin{eqnarray}
\eta &\stackrel{{\rm coev}\tens I}\longmapsto & 
\sum_\xi (\xi\bar\ra \tau(\langle \xi\rangle^L,\langle \xi\rangle)^{-1}\tens\hat\xi)
\tens \eta \cr
&\stackrel{{\Phi}}\longmapsto  & \sum_\xi \xi\bar\ra \tau(\langle \xi\rangle^L,\langle \xi\rangle)^{-1}
\tau(\langle \xi\rangle^L,\langle \eta\rangle)
\tens(\hat\xi \tens  \eta)  \cr  &\stackrel{I\tens{\rm eval}}\longmapsto  &
\sum_\xi \xi\bar\ra \tau(\langle \xi\rangle^L,\langle \xi\rangle)^{-1}
\tau(\langle \xi\rangle^L,\langle \eta\rangle)\ \delta_{\xi,\eta}\ =\ \eta\ . \cr
\hat\eta&\stackrel{I\tens{\rm coev}}\longmapsto & 
\sum_\xi \hat\eta\tens (\xi\bar\ra \tau(\langle \xi\rangle^L,\langle \xi\rangle)^{-1}\tens\hat\xi) \cr
&\stackrel{\Phi^{-1}}\longmapsto  & 
\sum_\xi (\hat\eta\bar\ra \tau(\langle \xi\rangle\ra 
\tau(\langle \xi\rangle^L,\langle \xi\rangle)^{-1}, \langle \xi\rangle^L)^{-1}
\tens \xi\bar\ra \tau(\langle \xi\rangle^L,\langle \xi\rangle)^{-1})\tens\hat\xi\ .\label{ill3}
\end{eqnarray}
Now we use the calculation
\begin{eqnarray*}
e &=& \langle \xi\rangle\tau(\langle \xi\rangle^L,\langle \xi\rangle)^{-1}\langle \xi\rangle^L\ =\ 
(\langle \xi\rangle \la \tau(\langle \xi\rangle^L,\langle \xi\rangle)^{-1})
(\langle \xi\rangle \ra \tau(\langle \xi\rangle^L,\langle \xi\rangle)^{-1}) \langle\xi\rangle^L\cr
&=& (\langle \xi\rangle \la \tau(\langle \xi\rangle^L,\langle \xi\rangle)^{-1})\ 
\tau(\langle \xi\rangle \ra \tau(\langle \xi\rangle^L,\langle \xi\rangle)^{-1},
\langle\xi\rangle^L)\ (
(\langle \xi\rangle \ra \tau(\langle \xi\rangle^L,\langle \xi\rangle)^{-1})\cdot \langle\xi\rangle^L)\ ,
\end{eqnarray*}
to rewrite the last line of (\ref{ill3}) as
\begin{eqnarray*}
&\stackrel{\Phi^{-1}}\longmapsto  & 
\sum_\xi (\hat\eta\bar\ra (\langle \xi\rangle \la \tau(\langle \xi\rangle^L,\langle \xi\rangle)^{-1})
\tens \xi\bar\ra \tau(\langle \xi\rangle^L,\langle \xi\rangle)^{-1})\tens\hat\xi \cr
&\stackrel{{\rm eval}\tens I}\longmapsto  & 
\sum_\xi \delta_{\xi,\eta}\tens\hat\xi\ =\ \hat\eta\ .
\end{eqnarray*}

\section{An algebra in the tensor category}
Take a group $X$ with subgroup $G$, and
a set of left coset representatives $M$ which contains $e$. 
We assume that $(M,\cdot)$ has the left division property, 
i.e.\  for all $t,s\in M$ there is a unique solution $p\in M$ to the equation $s\cdot p=t$.

We can combine the group action and the grading in the definition of ${\cal C}$
 by considering a single object $H$, 
a vector space spanned by a basis $\delta_s\tens u$ for $s\in M$ and $u\in G$. 
We suppose that $H$ is in the category ${\cal C}$, and define a map $\bar\ra:V\tens H\to V$
(for $V$ any object of ${\cal C}$) by
\[
\xi \bar\ra(\delta_s\tens u)\ =\ \delta_{s,\langle \xi\rangle }\ \xi\bar\ra u\ .
\]
If this map is to be a morphism in the category we must have $\langle \xi\rangle 
\cdot\langle \delta_s\tens u\rangle 
=\langle \xi\bar\ra u\rangle $ if $\langle \xi\rangle =s$, i.e.\
$s\cdot\langle \delta_s\tens u\rangle =s\ra u$. This can be solved uniquely for
$\langle \delta_s\tens u\rangle $ in $(M,\cdot)$ by left division.
 The action of $v\in G$ is given by (using $a=\langle \delta_s\tens u\rangle $)
\begin{eqnarray}
(\delta_s\tens u)\bar\ra v\ =\ \delta_{s\ra(a\la v)}\tens (a\la v)^{-1}uv\ .\label{etal}
\end{eqnarray}

\begin{propos} The action and grading on $H$ are consistent.
Further $\bar\ra:V\tens H\to V$ is a morphism in ${\cal C}$, for $V$ any object of ${\cal C}$. 
\end{propos}
\proof First we check that $\langle (\delta_s\tens u)\bar\ra v\rangle =
\langle \delta_s\tens u\rangle \ra v$. If we set $b=\langle (\delta_s\tens u)\bar\ra v\rangle $
then from (\ref{etal}), $(s\ra(a\la v))\cdot b=s\ra uv$,
where $a=\langle \delta_s\tens u\rangle $. If we apply $\ra v$ to the
equation $s\cdot\langle \delta_s\tens u\rangle =s\ra u$ and use the uniqueness
of the result of the left division process, we see $b=\langle \delta_s\tens u\rangle \ra v$. 

The grading on $H$ was defined so that $\bar\ra$ preserved the grades, so we only have to 
check the $G$-action. If we set $a=\langle \delta_s\tens u\rangle $ again, then
\[
(\xi\tens(\delta_s\tens u))\bar\ra v\ =\ \xi\bar\ra(a\la v)\tens 
(\delta_{s\ra(a\la v)}\tens (a\la v)^{-1}uv)\ ,
\]
and applying $\bar\ra$ to this gives
\[
\delta_{\langle \xi\rangle \ra(a\la v),s\ra(a\la v)}\ \xi\bar\ra uv\ ,
\]
which is just $(\xi\bar\ra (\delta_s\tens u))\bar\ra v$ as required.
\endproof

We would now like to give $H$ a multiplication so that $\bar\ra$ becomes an action of the algebra
$H$. Note that the result is not the usual semi-direct product multiplication. 

\begin{propos} The formula for the product $\mu$ for $H$ in ${\cal C}$ consistent with action
above is
\[
(\delta_s\tens u)(\delta_t\tens v)\ =\ \delta_{t,s\ra u}\ \delta_{s\ra\tau(a,b)}
\tens \tau(a,b)^{-1}uv\ ,
\]
where $a=\langle \delta_s\tens u\rangle $ and $b=\langle \delta_t\tens v\rangle $. 
\end{propos}
\proof We want the following equation to hold,
remembering to use $\Phi$ when we change the bracketing:
\begin{eqnarray}
(\xi\bar\ra(\delta_s\tens u))\bar\ra(\delta_t\tens v)\ =\ 
(\xi\bar\ra\tau(a,b))\bar\ra((\delta_s\tens u)(\delta_t\tens v))\ ,\label{etiv}
\end{eqnarray}
where $a=\langle \delta_s\tens u\rangle $ and $b=\langle \delta_t\tens v\rangle $. 
Now 
\begin{eqnarray*}
(\xi\bar\ra(\delta_s\tens u))\bar\ra(\delta_t\tens v) &=&
\delta_{s,\langle \xi\rangle }\ (\xi\bar\ra u)\bar\ra(\delta_t\tens v) \cr &=&
\delta_{s,\langle \xi\rangle }\ \delta_{t,\langle \xi\rangle \ra u}\ \xi\bar\ra uv\ ,
\end{eqnarray*}
and the two sides of (\ref{etiv}) agree by definition of the product above.
\endproof

\begin{propos} Multiplication $\mu:H\tens H\to H$ is a morphism in ${\cal C}$.
\end{propos}
\proof Set $\eta=\delta_s\tens u$, $\xi=\delta_t\tens v$, 
$a=\langle \eta\rangle $ and $b=\langle \xi\rangle $. 
For the grading, note that by definition $s\cdot a=s\ra u$, $t\cdot b=t\ra v$ and
$(s\ra\tau(a,b))\cdot\langle \eta\xi\rangle =s\ra uv$, under the assumption that $s\ra u=t$. 
But then $(s\cdot a)\cdot b=s\ra uv$, so 
$(s\ra\tau(a,b))\cdot(a\cdot b)=s\ra uv$ and we deduce that 
$\langle \eta\xi\rangle =a\cdot b$.

Now we check the action:
\[
(\eta\tens\xi)\bar\ra w\ =\ (\delta_{s\ra(a\la(b\la w))}\tens (a\la(b\la w))^{-1}u(b\la w))
\ \tens\ (\delta_{t\ra(b\la w)}\tens (b\la w)^{-1}vw)\ ,
\]
and multiplying these together gives
\[
\delta_{t,s\ra u}\ \ \delta_{s\ra\tau(a,b)(a\cdot b\la w)}\tens
(a\cdot b\la w)^{-1}\tau(a,b)^{-1}uvw\ ,
\]
which is the same as $(\eta\xi)\bar\ra w$. 
\endproof

\begin{propos} Multiplication $\mu:H\tens H\to H$ is associative in ${\cal C}$.
There is an identity $I$ for the multiplication and an algebra map $\eps:H\to k$
in the category, given by
\[
I\ =\ \sum_{t}\delta_t\tens e\ ,\quad \eps(\delta_s\tens u)\ =\ \delta_{s,e}\ .
\]
In terms of the action of $H$ on objects in ${\cal C}$, the identity 
$I$ has the trivial action on all objects. The action of $h\in H$ on the object $k$
is just multiplication by $\eps(h)$, and $\eps(I)=1$. 
\end{propos}
\proof Set $a=\langle \delta_s\tens u\rangle $, 
$b=\langle \delta_t\tens v\rangle $ and $c=\langle \delta_r\tens w\rangle $. Then
\begin{eqnarray*}
((\delta_s\tens u)(\delta_t\tens v))(\delta_r\tens w)&=&
\delta_{t,s\ra u}\ (\delta_{s\ra\tau(a,b)}
\tens \tau(a,b)^{-1}uv)(\delta_r\tens w) \cr
&=& \delta_{t,s\ra u} \delta_{r,s\ra uv}\ \delta_{s\ra\tau(a,b)\tau(a\cdot b,c)}\cr &&\quad \tens
\tau(a\cdot b,c)^{-1}\tau(a,b)^{-1}uvw\ , \cr
((\delta_s\tens u)\bar\ra\tau(b,c))((\delta_t\tens v)(\delta_r\tens w))&=&\delta_{r,t\ra v}\ 
(\delta_{s\ra(a\la\tau(b,c))}\tens (a\la\tau(b,c))^{-1}u\tau(b,c)) \cr &&\quad
(\delta_{t\ra\tau(b,c)}\tens \tau(b,c)^{-1}vw) \cr &=&
\delta_{t,s\ra u} \delta_{r,s\ra uv}\ 
\delta_{s\ra(a\la\tau(b,c))\tau(a\ra\tau(b,c),b\cdot c)}\cr &&\quad \tens
\tau(a\ra\tau(b,c),b\cdot c)^{-1}(a\la\tau(b,c))^{-1}uvw\ ,
\end{eqnarray*}
and these are equal by standard identities on $\tau$.

For the identity, note that $\langle I\rangle =e$, which is required as strictly the
identity is a morphism $:k\to H$ in the category. The rest is standard.
\endproof

\section{A braided tensor category}
Take a group $X$ with subgroup $G$, and
a set of left coset representatives $M$ which contains $e$. 
We consider a subcategory ${\cal D}$ of ${\cal C}$ with the additional structures of
a function $\bar\la:M\times V\to V$
and a $G$-grading, written $|\xi|\in G$ for $\xi$ in
 every object $V$ in the category ${\cal D}$. We require the following 
connections between the gradings and actions:
\begin{eqnarray}
|\eta\bar\ra u|\ =\ (\langle \eta\rangle \la u)^{-1}|\eta|u & , &
s\cdot\langle \eta\rangle \ =\ \langle s\bar\la\eta\rangle \cdot (s\ra|\eta|)\ ,\cr 
\tau(s,\langle \eta\rangle )^{-1}(s\la|\eta|) &=&
\tau(\langle s\bar\la\eta\rangle ,s\ra|\eta|)^{-1}
|s\bar\la\eta|\ .\label{connect}
\end{eqnarray}
The operation $\bar\la$ is an `action' of $M$, which we define to mean that
$ t \bar\la:V\to V$ is linear for all objects $V$ and all $t\in M$, and that
\begin{eqnarray}
 p\bar\la( t \bar\la\kappa) 
&=&
( p'\cdot t 
\bar\la\kappa)\bar\ra\tau(
p'\ra( t \la |\kappa|) 
, t \ra |\kappa| )^{-1}  \ ,\label{ten1}
\end{eqnarray}
where $p'=p\ra \tau(\langle  t \bar\la\kappa\rangle ,
 t \ra|\kappa|)\tau( t ,\langle \kappa\rangle )^{-1}$.
We also require a cross relation between the two actions,
\begin{eqnarray}
( s \bar\la\eta)\bar\ra(( s \ra|\eta|)\la u) &=&
( s \ra(\langle \eta\rangle \la u))\bar\la (\eta\bar\ra u)\ .\label{ten2}
\end{eqnarray}
The morphisms in the category ${\cal D}$ are linear maps preserving both gradings and both actions.

\begin{propos}\label{gra1} The connections between the gradings
and the actions are given by the following factorisations in $X$:
\begin{eqnarray*}
|s\bar\la\eta|^{-1}\langle s\bar\la\eta \rangle  &=& (s\ra|\eta|)
|\eta|^{-1}\langle \eta \rangle (s\ra|\eta|)^{-1}\ ,\cr
|\eta\bar\ra u|^{-1}\langle \eta\bar\ra u \rangle  &=& u^{-1}
|\eta|^{-1}\langle \eta \rangle u\ .
\end{eqnarray*}
\end{propos}
\proof Directly from the conditions above.
\endproof

Now we would like to make ${\cal D}$ into a tensor category. To do this we 
give the $G$-grading and action of $M$ on tensor products, and show that
the associator is a morphism. We define
\begin{eqnarray}
|\xi\tens\eta| &=& \tau(\langle \xi\rangle ,\langle \eta\rangle )^{-1}|\xi||\eta|\ ,\cr
 (s \ra\tau(\langle \eta\rangle ,\langle \kappa\rangle ) )
\bar\la (\eta\tens\kappa) &=&
( s \bar\la\eta)\bar\ra\tau( s \ra|\eta|,\langle \kappa\rangle ) 
\tau(\langle ( s \ra|\eta|)\bar\la \kappa\rangle , s \ra|\eta||\kappa|)^{-1}
\tens( s \ra|\eta|)\bar\la \kappa\ .\label{tend}
\end{eqnarray}

\begin{propos}\label{gra2} The gradings on the tensor product of objects $V\tens W$
 are given by the following factorisation
in $X$:
$
|\xi\tens\eta|^{-1}\langle \xi\tens\eta\rangle \ =\ |\eta|^{-1}|\xi|^{-1}
\langle \xi\rangle \ \langle \eta\rangle \ .
$
\end{propos}
\proof 
\[
|\xi\tens\eta|^{-1}\langle \xi\tens\eta\rangle \ =\ |\xi\tens\eta|^{-1}(\langle \xi\rangle 
\cdot\langle \eta\rangle )= |\xi\tens\eta|^{-1}\tau(\langle \xi\rangle ,\langle \eta\rangle )^{-1}
\langle \xi\rangle 
\langle \eta\rangle \ = |\eta|^{-1}|\xi|^{-1}
\langle \xi\rangle \ \langle \eta\rangle 
\]
\endproof

\begin{propos} The gradings on the tensor product are consistent with the actions,
as specified in (\ref{gra1}).
\end{propos}
\proof
First we check  the $G$-action. From (\ref{gra2}), for all $u\in G$,
\begin{eqnarray*}
|\xi\bar\ra(\langle \eta\rangle \la u)\tens \eta\bar\ra u|^{-1}
\langle \xi\bar\ra(\langle \eta\rangle \la u)\tens \eta\bar\ra u\rangle  &=&
|\eta\bar\ra u|^{-1}|\xi\bar\ra(\langle \eta\rangle \la u)|^{-1}
\langle \xi\bar\ra(\langle \eta\rangle \la u)\rangle  \ 
\langle \eta\bar\ra u\rangle   \cr &=&
|\eta\bar\ra u|^{-1}  (\langle \eta\rangle \la u)^{-1}   |\xi|^{-1}
\langle \xi\rangle  \  (\langle \eta\rangle \la u)\  
\langle \eta\bar\ra u\rangle   \cr &=&
u^{-1}|\eta|^{-1}|\xi|^{-1}
\langle \xi\rangle \ \langle \eta\rangle \  u\ .
\end{eqnarray*}
Now we check the $M$ action by considering the
grades of $(s \ra\tau(\langle \eta\rangle ,\langle \kappa\rangle ) )
\bar\la (\eta\tens\kappa)$. We set
 $u=\tau(s\ra|\eta|,\langle \kappa\rangle )\tau(\langle (s\ra|\eta|)
\bar\la\kappa\rangle ,s\ra|\eta||\kappa|)^{-1}$.
\begin{eqnarray*}
&& |(s\bar\la\eta)\bar\ra u\tens (s\ra|\eta|)\bar\la\kappa|^{-1}
\langle (s\bar\la\eta)\bar\ra u\tens (s\ra|\eta|)\bar\la\kappa\rangle \cr &&\quad =\ 
|(s\ra|\eta|)\bar\la\kappa|^{-1} |(s\bar\la\eta)\bar\ra u|^{-1}
\langle (s\bar\la\eta)\bar\ra u\rangle \  \langle (s\ra|\eta|)\bar\la\kappa\rangle 
\cr &&\quad =\ 
|(s\ra|\eta|)\bar\la\kappa|^{-1}u^{-1} |s\bar\la\eta|^{-1}
\langle s\bar\la\eta\rangle \, u\,  \langle (s\ra|\eta|)\bar\la\kappa\rangle 
\cr &&\quad =\ 
|(s\ra|\eta|)\bar\la\kappa|^{-1}u^{-1} (s\ra|\eta|) |\eta|^{-1}
\langle \eta\rangle \,  (s\ra|\eta|)^{-1}u\,  \langle (s\ra|\eta|)\bar\la\kappa\rangle \ .
\end{eqnarray*}
Now use the fact from(\ref{connect}) that $u|(s\ra|\eta|)\bar\la\kappa|=
(s\ra|\eta|)\la|\kappa|$ to get
\begin{eqnarray*}
&& |(s\bar\la\eta)\bar\ra u\tens (s\ra|\eta|)\bar\la\kappa|^{-1}
\langle (s\bar\la\eta)\bar\ra u\tens (s\ra|\eta|)\bar\la\kappa\rangle \cr &&\quad =\ 
((s\ra|\eta|)\la|\kappa|)^{-1} (s\ra|\eta|) |\eta|^{-1}
\langle \eta\rangle \,  (s\ra|\eta|)^{-1} ((s\ra|\eta|)\la|\kappa|) 
|(s\ra|\eta|)\bar\la\kappa|^{-1}
\,  \langle (s\ra|\eta|)\bar\la\kappa\rangle  \cr &&\quad =\ 
((s\ra|\eta|)\la|\kappa|)^{-1} (s\ra|\eta|) |\eta|^{-1}
\langle \eta\rangle \,  (s\ra|\eta|)^{-1} ((s\ra|\eta|)\la|\kappa|) 
(s\ra|\eta||\kappa|) |\kappa|^{-1}\langle \kappa\rangle (s\ra|\eta||\kappa|)^{-1}
\cr &&\quad =\ 
(s\ra|\eta||\kappa|) |\kappa|^{-1} |\eta|^{-1}
\langle \eta\rangle \,  \langle \kappa\rangle (s\ra|\eta||\kappa|)^{-1} \cr
 &&\quad =\ 
(s\ra\tau(\langle \eta\rangle ,\langle \kappa\rangle ) |\eta\tens\kappa|)\, |\kappa|^{-1} |\eta|^{-1}
\langle \eta\rangle \,  \langle \kappa\rangle \,
(s\ra\tau(\langle \eta\rangle ,\langle \kappa\rangle ) |\eta\tens\kappa|)^{-1}\ ,
\end{eqnarray*}
as required. \endproof

\begin{propos} The function $\bar \la$ applied to $V\tens W$ satisfies the condition
(\ref{ten1})
to be an $M$-action, and $\bar \la$ and $\bar \ra$ satisfy the cross relation
(\ref{ten2}) on $V\tens W$. 
\end{propos}
\proof Set $a=\langle \xi\rangle $ and $b=\langle \eta\rangle $, and begin with the cross relation,
with the following formula derived from the left hand side of (\ref{ten2}):
\begin{eqnarray}
&&((s\ra\tau(a,b))\bar\la(\xi\tens\eta))\bar\ra((s\ra\tau(a,b)|\xi\tens\eta|)\la u)\ =\ 
((s\bar\la \xi)\bar\ra v\tens (s\ra|\xi|)\bar\la\eta)\, \bar\ra ((s\ra|\xi||\eta|)\la u)  \cr
&&=\ (s\bar\la \xi)\bar\ra\, v(\langle (s\ra|\xi|)\bar\la\eta\rangle \la ((s\ra|\xi||\eta|)\la u))
\, \tens\, ((s\ra|\xi|)\bar\la\eta)\bar\ra ((s\ra|\xi||\eta|)\la u)  \cr
&&=\ (s\bar\la \xi)\bar\ra\, \tau(s\ra|\xi|,b)\ 
(\langle (s\ra|\xi|)\bar\la\eta\rangle \cdot(s\ra|\xi||\eta|)\la u)\ 
\tau(\langle (s\ra|\xi|)\bar\la\eta\rangle \ra ((s\ra|\xi||\eta|)\la u), s\ra|\xi||\eta|u)^{-1} \cr
&&\qquad \tens\, (s\ra|\xi|(b\la u))\bar\la(\eta\bar\ra u)\cr
&&=\ (s\bar\la \xi)\bar\ra\, \tau(s\ra|\xi|,b)\ 
((s\ra|\xi|)\cdot b\la u)\ 
\tau(\langle (s\ra|\xi|)\bar\la\eta\rangle \ra ((s\ra|\xi||\eta|)\la u), s\ra|\xi||\eta|u)^{-1} \cr
&&\qquad \tens\, (s\ra|\xi|(b\la u))\bar\la(\eta\bar\ra u) \ ,\label{ret2}
\end{eqnarray}
where $v=\tau(s\ra|\xi|,b)\, \tau(\langle (s\ra|\xi|)\bar\la\eta\rangle ,s\ra|\xi||\eta|)^{-1}$. 
This should be the same as the formula derived from the right hand side of (\ref{ten2}):
\begin{eqnarray}
(s\ra\tau(a,b)(a\cdot b\la u))\bar\la((\xi\tens\eta)\bar\ra u) &=&
 (s\ra\tau(a,b)(a\cdot b\la u))\bar\la(\xi\bar\ra(b\la u)\tens\eta\bar\ra u)  \cr
&=& (t\bar\la(\xi\bar\ra(b\la u)))\bar\ra w\tens (t\ra|\xi\bar\ra(b\la u)|)\bar\la
(\eta\bar\ra u)\ ,\label{ret}
\end{eqnarray}
where we have set
\begin{eqnarray*}
t &=& s\ra\tau(a,b)(a\cdot b\la u)\tau(a\ra(b\la u),b\ra u)^{-1}\ =\ s\ra(a\la(b\la u))\ ,\cr
w &=& \tau(t\ra|\xi\bar\ra(b\la u)|,b\ra u)\ \tau(
\langle  (t\ra|\xi\bar\ra(b\la u)|)\bar\la(\eta\bar\ra u) \rangle , t\ra |\xi\bar\ra(b\la u)|
|\eta\bar\ra u|)^{-1}\ .
\end{eqnarray*}
Now we can simplify some pieces of (\ref{ret}):
\begin{eqnarray*}
t\ra|\xi\bar\ra(b\la u)| &=& t\ra (a\la(b\la u))^{-1}|\xi|(b\la u)\ =\ 
s\ra |\xi|(b\la u) ,\cr
t\bar\la(\xi\bar\ra(b\la u)) &=& (s\ra(a\la(b\la u)))\bar\la(\xi\bar\ra(b\la u)) \ =\ 
 (s\bar\la\xi)\bar\ra ((s\ra|\xi|)\la(b\la u)) \cr
&=& (s\bar\la\xi)\bar\ra \tau(s\ra|\xi|,b)\  ((s\ra|\xi|)\cdot b\la u)\ 
\tau(s\ra |\xi|(b\la u),b\ra u)^{-1}\ ,\cr
w&=& \tau(s\ra |\xi|(b\la u),b\ra u)\, \tau(
\langle  (s\ra |\xi|(b\la u))\bar\la(\eta\bar\ra u) \rangle , s\ra |\xi|(b\la u)|\eta\bar\ra u|)^{-1}\cr
&=& \tau(s\ra |\xi|(b\la u),b\ra u)\, \tau(
\langle  ((s\ra |\xi|)\bar\la\eta)\bar\ra ((s\ra |\xi||\eta|)\la u) \rangle , 
s\ra |\xi||\eta|u)^{-1}\ .
\end{eqnarray*}
Substituting these in (\ref{ret}) gives the same result as (\ref{ret2}), as required.

Now we check the condition for the $M$-action. Begin with
\begin{eqnarray}
&&(p\ra \tau(\langle  s\bar\la\xi\rangle  \ra v, \langle  (s\ra|\xi|)\bar\la\eta\rangle  ))
\ \bar\la\ ((s\ra\tau(a,b))\bar\la(\xi\tens\eta)) \cr &&=\ 
(p\ra \tau(\langle  s\bar\la\xi\rangle  \ra v, \langle  (s\ra|\xi|)\bar\la\eta\rangle  ))
\ \bar\la\ ((s\bar\la\xi)\bar\ra v\tens(s\ra|\xi|)\bar\la\eta) \cr &&=\ 
(p\bar\la((s\bar\la \xi)\bar\ra v))\bar\ra z\ \tens\  (p\ra|(s\bar\la\xi)\bar\ra v|)\bar\la
((s\ra|\xi|)\bar\la\eta)\ ,\label{ret3}
\end{eqnarray}
where we have set
\[
z\ =\ \tau(p\ra|(s\bar\la\xi)\bar\ra v|, \langle  (s\ra|\xi|)\bar\la\eta\rangle )\ 
\tau(\langle   (p\ra|(s\bar\la\xi)\bar\ra v|)\bar\la
((s\ra|\xi|)\bar\la\eta)\rangle  , p\ra |(s\bar\la\xi)\bar\ra v||(s\ra|\xi|)\bar\la\eta|)^{-1}\ .
\]
We wish to show that (\ref{ret3}) is the same as
\begin{eqnarray}
&&(p'\cdot(s\ra\tau(a,b))\bar\la(\xi\tens\eta))\bar\ra \tau(p'\ra((s\ra\tau(a,b))\la |\xi\tens\eta|),
s\ra\tau(a,b)|\xi\tens\eta|)^{-1}  \cr
&&=\ 
(p'\cdot(s\ra\tau(a,b))\bar\la(\xi\tens\eta))\bar\ra \tau(p'\ra
(s\la\tau(a,b))^{-1}(s\la|\xi||\eta|),s\ra|\xi||\eta|)^{-1}
\ ,
\label{ret4}\end{eqnarray}
Where we have set
\begin{eqnarray*}
p' &=& p\ra \tau(\langle  s\bar\la\xi\rangle  \ra v, \langle  (s\ra|\xi|)\bar\la\eta\rangle  )\ 
\tau(\langle  (s\ra\tau(a,b))\bar\la(\xi\tens\eta)\rangle , s\ra\tau(a,b)|\xi\tens\eta|)\, 
\tau(s\ra\tau(a,b),a\cdot b)^{-1}  \cr
&=& p\ra ((\langle  s\bar\la\xi\rangle  \ra v)\la\tau(\langle  (s\ra|\xi|)\bar\la\eta\rangle ,
s\ra|\xi||\eta|))\cr &&\quad \tau(\langle  s\bar\la\xi\rangle  \ra v\tau(\langle  (s\ra|\xi|)\bar\la\eta\rangle ,
s\ra|\xi||\eta|), \langle  (s\ra|\xi|)\bar\la\eta\rangle \cdot s\ra|\xi||\eta|)\ 
\tau(s\ra\tau(a,b),a\cdot b)^{-1} \cr
&=& p\ra ((\langle  s\bar\la\xi\rangle  \ra v)\la v^{-1}\tau(s\ra|\xi|,b))
\ \tau(\langle  s\bar\la\xi\rangle  \ra \tau(s\ra|\xi|,b),  (s\ra|\xi|)\cdot b)\ 
\tau(s\ra\tau(a,b),a\cdot b)^{-1}\cr
&=& p\ra (\langle  s\bar\la\xi\rangle  \la v)^{-1}(\langle  s\bar\la\xi\rangle  \la 
\tau(s\ra|\xi|,b))
\  \tau(\langle  s\bar\la\xi\rangle  \ra \tau(s\ra|\xi|,b),  (s\ra|\xi|)\cdot b)\ 
\tau(s\ra\tau(a,b),a\cdot b)^{-1}\cr
&=& p\ra (\langle  s\bar\la\xi\rangle  \la v)^{-1}
\tau(\langle  s\bar\la\xi\rangle ,s\ra|\xi|)\  \tau(\langle  s\bar\la\xi\rangle \cdot (s\ra|\xi|),b)\ 
\tau(s\ra\tau(a,b),a\cdot b)^{-1}\cr
&=& p\ra (\langle  s\bar\la\xi\rangle  \la v)^{-1}\ 
\tau(\langle  s\bar\la\xi\rangle ,s\ra|\xi|)\  \tau(s\cdot a,b)\ 
\tau(s\ra\tau(a,b),a\cdot b)^{-1}\cr
&=& p\ra (\langle  s\bar\la\xi\rangle  \la v)^{-1}\ 
\tau(\langle  s\bar\la\xi\rangle ,s\ra|\xi|)\ \tau(s,a)^{-1}(s\la\tau(a,b))
\end{eqnarray*}
If we put $p''=p\ra (\langle  s\bar\la\xi\rangle  \la v)^{-1}\ 
\tau(\langle  s\bar\la\xi\rangle ,s\ra|\xi|)\, \tau(s,a)^{-1}$ then  (\ref{ret4}) becomes
\begin{eqnarray}
&&
(((p''\cdot s)\ra\tau(a,b))\bar\la(\xi\tens\eta))\ \bar\ra\  \tau(p''\ra
(s\la|\xi||\eta|),s\ra|\xi||\eta|)^{-1} \cr
&&=\ ((p''\cdot s\bar\la \xi)\bar\ra u'\tens ((p''\cdot s)\ra|\xi|)\bar \la\eta)
\ \bar\ra\  \tau(p''\ra
(s\la|\xi||\eta|),s\ra|\xi||\eta|)^{-1} \ ,
\label{ret5}\end{eqnarray}
where $u'=\tau((p''\cdot s)\ra|\xi|,b)\ \tau(\langle  ((p''\cdot s)\ra|\xi|)\bar \la\eta\rangle  ,
(p''\cdot s)\ra|\xi||\eta|)^{-1}$. 
We can simplify matters by
\begin{eqnarray*}
|(s\bar\la\xi)\bar\ra v| &=& (\langle  s\bar\la \xi\rangle \la v)^{-1}|s\bar\la \xi|v\ ,\cr
 (p\ra|(s\bar\la\xi)\bar\ra v|)\bar\la
((s\ra|\xi|)\bar\la\eta) &=& (q\cdot (s\ra|\xi|)\bar\la\eta)\bar\ra
\tau(q\ra ((s\ra|\xi|)\la|\eta|),s\ra|\xi||\eta|)^{-1}  \cr 
&=& (((p''\cdot s)\ra|\xi|)\bar\la\eta)\bar\ra
\tau(q\ra ((s\ra|\xi|)\la|\eta|),s\ra|\xi||\eta|)^{-1}\ , \cr
p\bar\la((s\bar\la \xi)\bar\ra v) &=& (c\bar\la(s\bar\la\xi))\bar\ra((c\ra|s\bar\la\xi|)\la v) \cr
&=& (p''\cdot s\bar\la\xi)\bar\ra \tau(p''\ra(s\la|\xi|),s\ra|\xi|)^{-1}
((c\ra|s\bar\la\xi|)\la v)\cr
&=& (p''\cdot s\bar\la\xi)\bar\ra \tau(q,s\ra|\xi|)^{-1}
(q\la v)\ ,
\end{eqnarray*}
where $c=p\ra(\langle  s\bar\la\xi\rangle \la v)^{-1}$ and
\begin{eqnarray*}
q &=& p\ra |(s\bar\la\xi)\bar\ra v| \, \tau(\langle  (s\ra|\xi|)\bar\la\eta\rangle ,s\ra|\xi||\eta|)
\tau(s\ra|\xi|,b)^{-1} \cr
&=& p\ra (\langle  s\bar\la\xi\rangle \la v)^{-1}|s\bar\la\xi|\ =\ p''\ra(s\la|\xi|)\ .
\end{eqnarray*}
Now we can rewrite (\ref{ret3}) as 
\[
 ((p''\cdot s)\bar\la\xi)\bar\ra \tau(q,s\ra|\xi|)^{-1}
(q\la v)z\ \tens\  ((p''\cdot s)\ra|\xi|)\bar\la\eta)\bar\ra
\tau(q\ra ((s\ra|\xi|)\la|\eta|),s\ra|\xi||\eta|)^{-1}\ ,
\]
so all we have to do now to show that this is equal to (\ref{ret5}) is to check that
\begin{eqnarray}
u'\, ( \langle  ((p''\cdot s)\ra|\xi|)\bar\la\eta\rangle  \la \tau(p''\ra
(s\la|\xi||\eta|),s\ra|\xi||\eta|)^{-1}) &=&\tau(q,s\ra|\xi|)^{-1}
(q\la v)z\ .\label{ret6}
\end{eqnarray}
To simplify what follows we shall use the substitutions
\begin{eqnarray*}
f\ =\ \langle  ((p''\cdot s)\ra|\xi|)\bar\la\eta\rangle  &,& n\ =\ \langle  (s\ra|\xi|)\bar\la\eta\rangle   \cr
g\ =\ p''\ra(s\la|\xi||\eta|) &,& h\ =\ s\ra|\xi||\eta|\ .
\end{eqnarray*}
If we use the result $p\ra|(s\bar\la \xi)\bar\ra v|=p''\ra(s\la|\xi|)v$, then we can rewrite
\begin{eqnarray*}
z &=& \tau(p''\ra(s\la|\xi|)v,n)\,
\tau(f \ra
\tau(g,h)^{-1},p''\ra(s\la|\xi|)v|(s\ra|\xi|)\bar\la\eta|)^{-1}\cr
&=& \tau(q\ra v,n)\,
\tau(f \ra
\tau(g,h)^{-1},g)^{-1}
\end{eqnarray*}
and (\ref{ret6}) becomes
\begin{eqnarray}
\tau((p''\cdot s)\ra|\xi|,b)\, \tau(f,g\cdot h)^{-1}(f\la\tau(g,h)^{-1}) &=&
\tau(q,s\ra|\xi|)^{-1}(q\la v)\, \tau(q\ra v,n)\, \tau(f \ra
\tau(g,h)^{-1},g)^{-1}   \cr
\tau(q,s\ra|\xi|)\, \tau(q\cdot(s\ra|\xi|),b) &=&
(q\la v)\, \tau(q\ra v,n)\, \tau((f\ra \tau(g,h)^{-1})\cdot g,h)
\ .\label{ret7}
\end{eqnarray}
Now we note some equations given by the grades:
\begin{eqnarray*}
((f\ra \tau(g,h)^{-1})\cdot g)\cdot h &=&
f\cdot(g\cdot h) \ =\  f\cdot((p''\cdot s)\ra|\xi||\eta|)\ =\ ((p''\cdot s)\ra|\xi|)\cdot b
\ =\ (q\cdot (s\ra|\xi|))\cdot b  \ , \cr
((q\ra v)\cdot n)\cdot h &=& (q\ra\tau(s\ra|\xi|,b))\cdot (n\cdot h)\ =\ 
(q\ra\tau(s\ra|\xi|,b))\cdot ((s\ra|\xi|)\cdot b)\ =\ (q\cdot (s\ra|\xi|))\cdot b\ .
\end{eqnarray*}
If we multiply (\ref{ret7}) on the right by $(q\cdot (s\ra|\xi|))\cdot b=
((f\ra \tau(g,h)^{-1})\cdot g)\cdot h$ we get
\begin{eqnarray}
\tau(q,s\ra|\xi|)\, (q\cdot(s\ra|\xi|))\, b &=&
(q\la v)\, \tau(q\ra v,n)\, ((f\ra \tau(g,h)^{-1})\cdot g)\, h
\ .\label{ret8}\end{eqnarray}
But we also have $(f\ra \tau(g,h)^{-1})\cdot g=(q\ra v)\cdot n$, so
\begin{eqnarray}
 q\,(s\ra|\xi|)\, b &=&
(q\la v)\, \tau(q\ra v,n)\, ((q\ra v)\cdot n)\, h\ =\ (q\la v)\, (q\ra v)\, n\, h  \cr
(s\ra|\xi|)\, b &=& v\, n\, h\ =\ v\, \tau(n,h)\, (n\cdot h)\ =\ \tau(s\ra|\xi|,b)\, ((s\ra|\xi|)\cdot b)
\ ,\label{ret9}\end{eqnarray}
which at last verifies (\ref{ret6}) and gives the answer!\endproof

\begin{theorem} When given the following structures, ${\cal D}$ is a braided tensor category:

The identity object is $k$, with trivial gradings and actions.

The associator $\Phi$ and the maps $l$ and $r$ are defined as for ${\cal C}$.

The braiding $\Psi:V\tens W \to W\tens V$ is defined by
$\Psi(\xi\tens\eta)=\langle \xi\rangle \bar\la\eta\tens\xi\bar\ra|\eta|$.
\end{theorem}
\proof The following lemmas.\endproof

\begin{lemma} The associator $\Phi$ is a morphism in the category ${\cal D}$.
\end{lemma}
\proof We begin by checking the $G$-grade.
\begin{eqnarray*}
|(\xi\tens\eta)\tens\kappa| &=& 
\tau(\langle \xi\rangle \cdot  \langle \eta\rangle ,\langle \kappa\rangle )^{-1}|\xi\tens\eta||\kappa|
\cr   &=& 
\tau(\langle \xi\rangle \cdot  \langle \eta\rangle ,\langle \kappa\rangle )^{-1}
\tau(\langle \xi\rangle , \langle \eta\rangle )^{-1}
|\xi||\eta||\kappa|  \cr
|\xi\bar\ra \tau(\langle \eta\rangle , \langle \kappa\rangle )\tens(\eta\tens\kappa)| &=& 
\tau(\langle \xi\rangle \ra \tau(\langle \eta\rangle , \langle \kappa\rangle ),
\langle \eta\rangle \! \cdot\! \langle \kappa\rangle )^{-1} 
|\xi\bar\ra \tau(\langle \eta\rangle , \langle \kappa\rangle )||\eta\tens\kappa|  \cr
&=& 
\tau(\langle \xi\rangle \ra \tau(\langle \eta\rangle , \langle \kappa\rangle ),
\langle \eta\rangle \! \cdot\! \langle \kappa\rangle )^{-1} 
(\langle \xi\rangle \la \tau(\langle \eta\rangle , \langle \kappa\rangle )  )^{-1}
|\xi||\eta||\kappa| \ .
\end{eqnarray*}
These are equal by the properties of $\tau$.

Now we check the $M$-action.  Set $\langle \xi\rangle =a$, 
$\langle \eta\rangle =b$ and $\langle \kappa\rangle =c$. To begin,
\begin{eqnarray}
(s\ra\tau(a\cdot b,c))\bar\la((\xi\tens\eta)\tens\kappa) &=&
(s\bar\la (\xi\tens\eta))\bar\ra w\tens (s\ra|\xi\tens\eta|)\bar\la \kappa   \cr
 &=&((t\bar\la\xi)\bar\ra u\tens(t\ra|\xi|)\bar\la \eta)\bar\ra w
\tens (s\ra|\xi\tens\eta|)\bar\la \kappa   \cr
 &=&((t\bar\la\xi)\bar\ra u (\langle  (t\ra|\xi|)\bar\la \eta \rangle \la w) \tens
((t\ra|\xi|)\bar\la \eta)\bar\ra w)\tens (s\ra|\xi\tens\eta|)\bar\la \kappa    ,\label{ll1}
\end{eqnarray}
where $s=t\ra\tau(a,b)$ and 
\begin{eqnarray*}
w &=& \tau(s\ra|\xi\tens\eta|,c)\, \tau(\langle  (s\ra|\xi\tens\eta|)\bar\la \kappa \rangle ,
s\ra|\xi\tens\eta||\kappa|)^{-1}\ ,\cr
u &=& \tau(t\ra|\xi|,b)\, \tau(\langle  (t\ra|\xi|)\bar\la \eta \rangle ,t\ra|\xi||\eta|)^{-1}\ .
\end{eqnarray*}
If we set  $\xi'=\xi\bar\ra\tau(b,c)$ and $p\ra\tau(a\ra\tau(b,c),b\cdot c)=s\ra\tau(a\cdot b,c)$,
we would like (\ref{ll1}) to equal 
\begin{eqnarray}
\Phi^{-1}((p\ra\tau(a\ra\tau(b,c),b\cdot c))\bar\la(\xi'\tens(\eta\tens\kappa))) &=& \Phi^{-1}(
(p\bar\la \xi')\bar\ra v\tens (p\ra|\xi'|)\bar\la(\eta\tens\kappa))  \cr &=&
\Phi^{-1}(
(p\bar\la \xi')\bar\ra v\tens (  (q\bar\la\eta)\bar\ra z\tens (q\ra|\eta|)\bar\la\kappa)) \cr &=&
((p\bar\la \xi')\bar\ra vx\tens   (q\bar\la\eta)\bar\ra z)\tens (q\ra|\eta|)\bar\la\kappa\ ,\label{ll2}
\end{eqnarray}
where $q\ra\tau(b,c)=p\ra|\xi'|$ and
\begin{eqnarray*}
v &=& \tau(p\ra|\xi'|,b\cdot c)\, \tau(\langle  (p\ra|\xi'|)\bar\la(\eta\tens\kappa) \rangle ,
p\ra|\xi'||\eta\tens\kappa|  )^{-1}\ ,\cr
z&=& \tau(q\ra|\eta|,c)\, \tau( \langle  (q\ra|\eta|)\bar\la\kappa  \rangle ,
q\ra|\eta||\kappa|)^{-1}\ , \cr
x&=& \tau(\langle  q\bar\la\eta  \rangle \ra z,
\langle  (q\ra|\eta|)\bar\la\kappa  \rangle )^{-1}\ .
\end{eqnarray*}
Using the usual identities on $\tau$,
\[
p\ra|\xi'|\ =\ p\ra(a\la\tau(b,c))^{-1}|\xi|\tau(b,c)\ =\ t\ra|\xi|\tau(b,c)\  ,
\]
or equivalently $q=t\ra|\xi|$. Then $s\ra|\xi\tens\eta|=q\ra|\eta|$, so the third terms
in (\ref{ll1}) and (\ref{ll2})
are equal. Also we get $z=w$, so $((t\ra|\xi|)\bar\la \eta)\bar\ra w=(q\bar\la\eta)\bar\ra z$, 
and the second terms are equal. Next by using the cross relation on the first term of (\ref{ll2}),
\[
(p\bar\la \xi')\bar\ra vx\ =\ (p\bar\la (\xi\bar\ra\tau(b,c)))\bar\ra vx
\ =\ 
(t\bar\la \xi)\bar\ra (q\la\tau(b,c))vx
\]
Now we are left with the task of showing that
\begin{eqnarray}
(q\la\tau(b,c))vx &=& u(\langle  q\bar\la\eta  \rangle \la z)\ .\label{ll3}
\end{eqnarray}
If we use the formula
$
\langle  (q\ra\tau(b,c))\bar\la(\eta\tens\kappa)  \rangle \, =\, 
(\langle  q\bar\la\eta  \rangle \ra z)\cdot
\langle ( q\ra|\eta|)\bar\la\kappa  \rangle 
$
then
\begin{eqnarray*}
vx &=& \tau(q\ra\tau(b,c),b\cdot c)\  \tau((\langle  q\bar\la\eta  \rangle \ra z)\cdot
\langle ( q\ra|\eta|)\bar\la\kappa  \rangle ,q\ra|\eta||\kappa|)^{-1}\ 
\tau(\langle  q\bar\la\eta  \rangle \ra z,
\langle  (q\ra|\eta|)\bar\la\kappa  \rangle )^{-1}  \cr
&=& \tau(q\ra\tau(b,c),b\cdot c)\  \tau(
\langle  q\bar\la\eta  \rangle \ra z\, \tau(\langle ( q\ra|\eta|)\bar\la\kappa  \rangle ,
q\ra|\eta||\kappa|)
,(q\ra|\eta|)\cdot c)^{-1}  \cr &&\quad ((\langle  q\bar\la\eta  \rangle \ra z)\la
\tau(\langle ( q\ra|\eta|)\bar\la\kappa  \rangle ,
q\ra|\eta||\kappa|))^{-1}
\end{eqnarray*}
where we have used $\langle ( q\ra|\eta|)\bar\la\kappa  \rangle \cdot (q\ra|\eta||\kappa|)=
(q\ra|\eta|)\cdot c$. Then
\begin{eqnarray*}
(q\la\tau(b,c))vx &=& (q\la\tau(b,c))\  \tau(q\ra\tau(b,c),b\cdot c)\  \tau(
\langle  q\bar\la\eta  \rangle \ra \tau(q\ra|\eta|,c)
,(q\ra|\eta|)\cdot c)^{-1}  \cr &&\quad ((\langle  q\bar\la\eta  \rangle \ra z)\la
\tau(\langle ( q\ra|\eta|)\bar\la\kappa  \rangle ,
q\ra|\eta||\kappa|))^{-1} \cr  &=& \tau(q,b)\  \tau(q\cdot b,c)\  \tau(
\langle  q\bar\la\eta  \rangle \ra \tau(q\ra|\eta|,c)
,(q\ra|\eta|)\cdot c)^{-1}  \cr &&\quad ((\langle  q\bar\la\eta  \rangle \ra z)\la
\tau(\langle ( q\ra|\eta|)\bar\la\kappa  \rangle ,
q\ra|\eta||\kappa|))^{-1} \cr
 &=& 
\tau(q,b)\  \tau(\langle  q\bar\la\eta  \rangle ,q\ra|\eta|)^{-1}\ (\langle  q\bar\la\eta  \rangle 
\la \tau(q\ra|\eta|,c))
  \cr &&\quad ((\langle  q\bar\la\eta  \rangle \ra z)\la
\tau(\langle ( q\ra|\eta|)\bar\la\kappa  \rangle ,
q\ra|\eta||\kappa|))^{-1}  \cr
 &=& 
u\ (\langle  q\bar\la\eta  \rangle 
\la \tau(q\ra|\eta|,c))\  ((\langle  q\bar\la\eta  \rangle 
\ra \tau(q\ra|\eta|,c))\la \tau(\langle ( q\ra|\eta|)\bar\la\kappa  \rangle ,
q\ra|\eta||\kappa|)^{-1}  )
\cr &=&
u (\langle  q\bar\la \eta \rangle \la z)\ ,
\end{eqnarray*}
as required, where we have used $\langle  q\bar\la\eta  \rangle \cdot 
(q\ra|\eta|)=q\cdot b$. 
\endproof

\begin{lemma} The maps $l_V$ and $r_V$ are morphisms in the category ${\cal D}$.
\end{lemma}
\proof This is reasonably simple from the definitions, rembering that $\tau(e,s)=\tau(s,e)=e$
for all $s\in M$. Only the $G$-grade and the $M$-action need be checked.\endproof

\begin{lemma} The map $\Psi:V\tens W \to W\tens V$ defined by
$\Psi(\xi\tens\eta)=\langle \xi\rangle \bar\la\eta\tens\xi\bar\ra|\eta|$
is a morphism in the category.
\end{lemma}
\proof First we check the grades, using (\ref{gra1}):
\begin{eqnarray*}
|\xi\bar\ra|\eta||^{-1}|\langle \xi\rangle \bar\la\eta|^{-1}
 \langle \langle \xi\rangle \bar\la\eta\rangle \  \langle \xi\bar\ra|\eta|\rangle  &=&
|\xi\bar\ra|\eta||^{-1}    (\langle \xi\rangle \ra|\eta|)        |\eta|^{-1}
 \langle \eta\rangle \  (\langle \xi\rangle \ra|\eta|)  ^{-1} \langle \xi\bar\ra|\eta|\rangle  \cr
&=&
|\xi\bar\ra|\eta||^{-1}    \langle \xi\bar\ra|\eta|\rangle        |\eta|^{-1}
 \langle \eta\rangle \ \cr
&=&
|\eta|^{-1}\ |\xi|^{-1}    \langle \xi\rangle    |\eta|\,    |\eta|^{-1}
 \langle \eta\rangle \ \ =\ |\eta|^{-1}\ |\xi|^{-1}    \langle \xi\rangle    \,
 \langle \eta\rangle \ .
\end{eqnarray*}

Now we check the $G$ action:
\begin{eqnarray*}
(\Psi(\xi\tens\eta))\bar\ra u &=&
(\langle \xi\rangle \bar\la\eta)\bar\ra((\langle \xi\rangle \ra|\eta|)\la u)\tens\xi\bar\ra|\eta|u\ ,
\cr 
\Psi((\xi\tens\eta)\bar\ra u)&=& \Psi(\xi\bar\ra(\langle \eta\rangle \la u)\tens \eta\bar\ra u)\ =\ 
(\langle \xi\rangle \ra(\langle \eta\rangle \la u))\bar\la (\eta\bar\ra u) \tens 
\xi\bar\ra(\langle \eta\rangle \la u)|\eta\bar\ra u|\ .
\end{eqnarray*}
The first terms are equal by the cross relation (\ref{ten2}), and the second terms are equal by the 
connections
between the grades and the relations (\ref{connect}).

Now we check the $M$ action. Set $a=\langle \xi\rangle $ and $b=\langle \eta\rangle $. 
\begin{eqnarray}
\Psi((s\ra\tau(a,b))\bar\la (\xi\tens\eta)) &=&
\Psi((s\bar\la \xi)\bar\ra v\tens (s\ra|\xi|)\bar\la \eta)\cr &=&
(\langle s\bar\la \xi\rangle \ra v)\bar\la((s\ra|\xi|)\bar\la \eta)\tens
(s\bar\la \xi)\bar\ra v|(s\ra|\xi|)\bar\la \eta|\ ,\label{ban1}
\end{eqnarray}
where $v=\tau(s\ra|\xi|,b)\tau(\langle (s\ra|\xi|)\bar\la \eta\rangle ,s\ra|\xi||\eta|
)^{-1}$. Note that $v|(s\ra|\xi|)\bar\la \eta|=(s\ra|\xi|)\la |\eta|$. Then
(\ref{ban1}) should be the same as
\begin{eqnarray}
(s\ra\tau(a,b))\bar\la\Psi (\xi\tens\eta) &=&
(s\ra\tau(a,b))\bar\la(a\bar\la \eta\tens\xi\bar\ra |\eta|)  \cr
&=& (p\bar\la(a\bar\la\eta))\bar\ra w\tens (p\ra|a\bar\la\eta|)\bar\la(\xi\bar\ra |\eta|)\ ,
\label{ban2}\end{eqnarray}
where $s\ra\tau(a,b)=p\ra\tau(\langle a\bar\la \eta\rangle ,a\ra|\eta|)$ and
\[
w\ =\ \tau(p\ra|a\bar\la\eta|,a\ra|\eta|)\, \tau(\langle 
(p\ra|a\bar\la\eta|)\bar\la(\xi\bar\ra |\eta|) \rangle , p\ra|a\bar\la\eta||\xi\bar\ra |\eta||)^{-1}\ .
\]
Then
$p\ra|a\bar\la \eta|=s\ra(a\la|\eta|)$, so the second term of (\ref{ban1})
can be written $(s\bar\la \xi)\bar\ra ((s\ra|\xi|)\la |\eta|)$, and the
second term of (\ref{ban2}) as $(s\ra(a\la|\eta|))\bar\la(\xi\bar\ra |\eta|)$.
These are equal
by the cross relation (\ref{ten2}). For the first terms,
\begin{eqnarray*}
(p\bar\la(a\bar\la \eta))\bar\ra w &=&
 ((s\cdot a)\bar \la\eta)\bar\ra\tau(s\ra(a\la|\eta|),a\ra|\eta|)^{-1}w \cr
&=&  ((s\cdot a)\bar \la\eta)\bar\ra \tau( \langle s\bar\la \xi\rangle \ra ((s\ra|\xi|)\la|\eta|),
s\ra(a\la|\xi|)|\xi\bar\ra|\eta||)^{-1}\ ,
\cr
(\langle s\bar\la \xi\rangle \ra v)\bar\la((s\ra|\xi|)\bar\la \eta) &=&
((q\cdot (s\ra|\xi|))\bar\la \eta)\bar\ra \tau( q\ra((s\ra|\xi|)\la|\eta|),s\ra|\xi||\eta|)^{-1} \ ,
\end{eqnarray*}
where \[
q\ =\ \langle s\bar\la \xi\rangle \ra\,  v\ \tau(\langle  (s\ra|\xi|)\bar\la \eta \rangle ,
s\ra|\xi||\eta|)\ \tau(s\ra|\xi|,b)^{-1}\ =\ \langle s\bar\la \xi\rangle \ .
\]
Now we see that $q\cdot (s\ra|\xi|)=s\cdot a$. 
\endproof

\begin{lemma} The map $\Psi$ satisfies the hexagon identities.
\end{lemma}
\proof Set $\langle \xi\rangle =a$, $\langle \eta\rangle =b$ and $\langle \kappa\rangle =c$.
The following two compositions can be seen to be equal:
\begin{eqnarray}
(\xi\tens\eta)\tens\kappa &\stackrel{\Phi}\longmapsto& \xi\bar\ra\tau(b,c)
\tens(\eta\tens\kappa)  \cr
&\stackrel{\Psi}\longmapsto& (a \ra\tau(b,c) )
\bar\la (\eta\tens\kappa)  \tens 
  \xi\bar\ra\tau(b,c)   | \eta\tens\kappa  |\cr
&=& ((a\bar\la\eta)\bar\ra u\tens (a\ra|\eta|)\bar\la\kappa)\tens \xi\bar\ra|\eta||\kappa|\cr
&\stackrel{\Phi}\longmapsto& 
(a\bar\la\eta)\bar\ra u\,\tau(\langle (a\ra|\eta|)\bar\la\kappa\rangle ,a\ra|\eta||\kappa)
\tens ((a\ra|\eta|)\bar\la\kappa\tens \xi\bar\ra|\eta||\kappa|)\ ,\cr 
(\xi\tens\eta)\tens\kappa &\stackrel{\Psi\tens I}\longmapsto&
 (a\bar\la\eta\tens\xi\bar\ra|\eta|)\tens\kappa
 \cr &\stackrel{\Phi}\longmapsto& 
(a\bar\la\eta)\bar\ra\tau(a\ra|\eta|,c) 
\tens(\xi\bar\ra|\eta|  \tens \kappa)
 \cr &\stackrel{I\tens\Psi}\longmapsto& 
(a\bar\la\eta)\bar\ra\tau(a\ra|\eta|,c) 
\tens((a\ra|\eta|)\bar\la \kappa  \tens \xi\bar\ra|\eta||\kappa|)\ ,\label{br1}
\end{eqnarray}
where $u=\tau(a\ra|\eta|,c)\,
\tau(\langle (a\ra|\eta|)\bar\la\kappa\rangle ,a\ra|\eta||\kappa|)^{-1}$.

The hexagon identity for the inverse associator asserts that the following should
be equal:
\begin{eqnarray}
\xi\tens(\eta\tens\kappa) &\stackrel{I\tens\Psi}\longmapsto& 
\xi\tens(b\bar\la\kappa\tens\eta\bar\ra|\kappa|)
\cr &\stackrel{\Phi^{-1}}\longmapsto& 
(  \xi\bar\ra\tau(\langle b\bar\la\kappa\rangle ,b\ra|\kappa|)^{-1}
\tens b\bar\la\kappa)\tens\eta\bar\ra|\kappa|
\cr &\stackrel{\Psi\tens I}\longmapsto& 
(   (a\ra \tau(\langle b\bar\la\kappa\rangle ,
b\ra|\kappa|)^{-1})\bar\la(b\bar\la\kappa) \tens
\xi\bar\ra \tau(\langle b\bar\la\kappa\rangle ,b\ra|\kappa|)^{-1}
|b\bar\la\kappa|)\tens \eta\bar\ra|\kappa|\ ,
\cr
\xi\tens(\eta\tens\kappa) &\stackrel{\Phi^{-1}}\longmapsto& 
(\xi\bar\ra\tau(b,c)^{-1}\tens\eta)\tens\kappa
\cr &\stackrel{\Psi}\longmapsto& 
  ((a\ra \tau(b,c)^{-1})\cdot b)
\bar\la\kappa \tens(
\xi\bar\ra\tau(b,c)^{-1}(b\la |\kappa|) 
\tens \eta\bar\ra|\kappa| )
\cr &\stackrel{\Phi^{-1}}\longmapsto& 
(( ((a\ra \tau(b,c)^{-1})\cdot b)
\bar\la\kappa)\bar\ra\tau(
a\ra\tau(b,c)^{-1}(b\la |\kappa|) 
,b\ra |\kappa| )^{-1}      \cr && \quad\quad \tens
\xi\bar\ra\tau(b,c)^{-1}(b\la |\kappa|) )
\tens \eta\bar\ra|\kappa| \ .\label{br2}
\end{eqnarray}
The third terms in (\ref{br2}) are equal, and the second terms can be seen to be equal by
the formula from (\ref{connect}),
$
\tau( b ,\langle \kappa\rangle )^{-1}( b \la |\kappa|) =
\tau(\langle  b \bar\la\kappa\rangle , b\ra|\kappa|)^{-1}
| b \bar\la\kappa|
$.
For the first terms use the condition for an $M$ action (\ref{ten1}) to get
\[
  ( a \ra \tau(\langle  b \bar\la\kappa\rangle ,
 b \ra|\kappa|)^{-1})\bar\la( b \bar\la\kappa) 
\ =\ 
( ( ( a \ra \tau( b ,c)^{-1})\cdot b )
\bar\la\kappa)\bar\ra\tau(
 a \ra\tau( b ,c)^{-1}( b \la |\kappa|) 
, b \ra |\kappa| )^{-1}  \ ,
\]
as required.\endproof

\begin{propos}\label{nat2} The braiding is a natural transformation
between the tensor product and its opposite in ${\cal D}$.
\end{propos} 
\proof The statement just means that the following diagram commutes for
all morphisms $\theta:V\to \tilde V$ and $\phi:W\to \tilde W$:
\[\begin{array}{ccc}
 V\tens W & \stackrel{\Psi_{VW}}\longrightarrow & W\tens V \\
\downarrow\ \theta\tens\phi & & \downarrow\ \phi\tens\theta \\
\tilde V\tens \tilde W & \stackrel{\Psi_{\tilde V\tilde W}}\longrightarrow & \tilde W\tens \tilde V
\end{array}\ ,\]
 This is simple to check, remembering that the morphisms preserve the grades
and actions. \endproof

\section{A double construction}
Take a group $X$ with subgroup $G$, and
a set of left coset representatives $M$.

\begin{defin}\label{tyy} The set $Y$, which is identical to $X$, is given 
a binary operation $\circ$ defined by
\[
(us)\circ(vt)\ =\ vust\ =\ vu\,\tau(s,t)(s\cdot t)\quad u,v\in G\ ,\ s,t\in M\ .
\]
Define the functions $\tilde\ra:Y\times X\to Y$ and $\tilde\tau:Y\times Y\to X$
 by $y\tilde\ra x=x^{-1}yx$ and $\tilde\tau(vt,wp)=\tau(t,p)$. The function
 $\tilde\la:Y\times X\to X$ is defined by
\[
vt\tilde\la wp\ =\ v^{-1}wpv'\ =\ twp{t'}^{-1}\ ,\quad{\rm where}\quad
vt\tilde\ra wp\ =\ v't'\ \quad v'\in G\ , \ t'\in M\ .
\]
\end{defin} 

\begin{propos} The maps $\tilde\ra$, $\tilde\la$ and $\tilde\tau$ satisfy
the six conditions listed in (\ref{rels}), with $(Y,\circ)$ taking the place of $(M,\cdot)$,
and the group $X$ taking the place of $G$.
\end{propos} 
\proof The fourth is immediate. For the sixth, consider
\[
(\tau(t,p)^{-1}us\,\tau(t,p))\circ((vt)\circ(wp)) \ =\  wvus\,\tau(t,p)(t\cdot p)\ =\ 
wvustp\ =\ ((us)\circ(vt))\circ(wp)
\ .
\]
For the second condition we start with
\[
v'u's't'\ =\ (u's')\circ(v't')\ =\ ((us)\circ(vt))\tilde\ra wp\ =\ p^{-1}w^{-1}vustwp\ ,
\]
where $(vt)\tilde\ra wp\ =\ p^{-1}w^{-1}vtwp\ =\ v't'$. 
From this we deduce that 
\[
u's'\ =\ {v'}^{-1}p^{-1}w^{-1}vustwp{t'}^{-1}\ =\ us\tilde\ra v^{-1}wpv'\ ,
\]
as  $v^{-1}wpv'=twp{t'}^{-1}$. 

For the third condition, we have $vt\tilde\la wpus=v^{-1}wpusv''$, 
where
\[
vt\tilde\ra wpus\ =\ v''t''\ =\ v't'\tilde\ra us\ \quad v''\in G\ , \ t''\in M\ .
\]
Then  $vt\tilde\la wpus=v^{-1}wp{v'}{v'}^{-1}usv''=(vt\tilde\la wp)((vt\tilde\ra wp)\tilde\la us)$. 

For the fifth condition, begin with $wp\tilde\ra\tilde\tau(us,vt)=\tau(s,t)^{-1}wp\tau(s,t)$,
so 
\begin{eqnarray*}
\tilde\tau( wp\tilde\ra\tilde\tau(us,vt), us\circ vt) &=& \tau(p\ra \tau(s,t), s\cdot t)\ ,
\cr
wp\tilde\la\tau(s,t)=p\,\tau(s,t)(p\ra \tau(s,t))^{-1} &=& p\la \tau(s,t)\ ,
\end{eqnarray*}
and from these we verify the fifth condition, which is
\[
(wp\tilde\la\tilde\tau(us,vt))\  \tilde\tau( wp\tilde\ra\tilde\tau(us,vt), us\circ vt)\ =\ 
\tilde\tau(wp,us)\  \tilde\tau((wp)\circ(us),vt)\ .
\]

For the first condition, begin with
\[
us\tilde\la (vt\tilde\la wp)\ =\ us\tilde\la v^{-1}wpv'\ =\ u^{-1}v^{-1}wpv'u'\ ,
\]
where $us\tilde\ra v^{-1}wpv'=us\tilde\ra(vt\tilde\la wp)=u's'$. Next
\[
(us\circ vt)\tilde\la wp\ =\ vu\tau(s,t)(s\cdot t)\tilde\la wp\ =\ 
\tau(s,t)^{-1}u^{-1}v^{-1}wp u''\ ,
\]
where, by the second condition,
 $(us\circ vt)\tilde\ra wp=u''s''=(us\tilde\ra(vt\tilde\la wp))\circ(vt\tilde\ra wp)=
u's'\circ v't'$, so we get $u''=v'u'\tau(s',t')$,
and hence verify the first condition. 
\endproof

\begin{propos} The element $e_y=f_m^{-1}e_m=e$ is a left identity for $Y$ (note it is not in
general a right identity), and the operation $(Y,\circ)$ has the right division property.
The corresponding left inverse is given by the formula 
$(vt)^L=v^{-1}t^{-1}$
(for $v\in G$ and $t\in M$).
\end{propos} 
\proof  To show that $e\in Y$ is a left inverse, note that $e\circ us=ues=us$ for all 
$us\in Y$ ($u\in G$ and $s\in M$). 

To check right division we have to check that there is a unique solution $wp\in Y$
($w\in G$ and $p\in M$) to the equation $wp\circ us=vt$.  The equation gives
$uw\,\tau(p,s)(p\cdot s)=vt$, and we can solve the equation $p\cdot s=t$
to give a unique value of $p$. Now $w=u^{-1}v\,\tau(p,s)^{-1}$. 

To check the formula for the left identity, $(v^{-1}t^{-1})\circ vt=vv^{-1}t^{-1}t=e=e_y$. 
\endproof

\begin{propos} If we define $f_y=e_m=f_m\in X$, we see that the conditions in 
 (\ref{inv}) are satisfied,
using $f_y$, $e_y$, $X$ and $(Y,\circ)$ instead of $f_m$, $e_m$, $G$ and $(M,\cdot)$. 
\end{propos} 
\proof  For the first condition, note that $e_y\tilde\ra x=x^{-1}ex=e_y$. This implies the second
condition, $e_y\tilde\la x=e_mxe_m^{-1}$. 

For the fourth condition, $us\tilde\ra e=us$, and then $us\tilde\la e=ses^{-1}=e$, 
the third condition. 

For the fifth condition, $\tilde\tau(e_y,us)=\tau(e_m,s)=f_m=f_y$. 

For the sixth condition, $us\tilde\ra f_y^{-1}=f_musf_m^{-1}=f_mu(s\la f_m^{-1})(s\ra f_m^{-1})$,
so $\tilde\tau(us\tilde\ra f_y^{-1},e_y)=\tau(s\ra f_m^{-1},e_m)=(s\la f_m^{-1})^{-1}$. 
Then $us\tilde\la f_y^{-1}=sf_m^{-1}(s\ra f_m^{-1})^{-1}=s\la f_m^{-1}$, as required. 

For the seventh condition, $(us\tilde\ra f_y^{-1})\circ e_m=f_m^{-1}f_musf_m^{-1}e_m=us$. \endproof

Now we return to the case where $e\in M$ for simplicity. 
We introduce a $Y$ valued grading on the objects of ${\cal D}$ by $\|\xi\|=|\xi|^{-1}\langle \xi\rangle $.
From our previous calculations we know that $\|\eta\bar\ra u\|=\|\eta\|\tilde\ra u$,
$\|s\bar\la\eta\|=\|\eta\|\tilde\ra(s\ra|\eta|)^{-1}$ and
 $\|\xi\tens\eta\|=\|\xi\|\circ\|\eta\|$. 

\begin{propos} The map $\hat\ra:V\times X\to V$ defined by $\xi\hat\ra us=(\xi\bar\ra u)\hat\ra s$
 ($u\in G$ and $s\in M$), where
\[
\xi\hat\ra s\ =\ ((s^L\ra|\xi|^{-1})\bar\la\xi)\bar\ra \tau(s^L,s)\ ,
\]
is a right action of the group $X$ on $V$, 
where $V$ is any object in ${\cal D}$. Further 
$\|\xi\hat\ra us\|=\|\xi\|\tilde\ra us$. 
\end{propos} 
\proof First consider the grading;
\[
\|\xi\hat\ra s\|\ =\ \| (s^L\ra|\xi|^{-1})\bar\la\xi \|\tilde\ra \tau(s^L,s)\ =\ 
\|\xi\|\tilde\ra s^{L-1} \tau(s^L,s)\ =\ \|\xi\|\tilde\ra s\ ,
\]
since $s^Ls=\tau(s^L,s)$. 

Now we wish to show that
\[
(\xi\hat\ra us)\hat\ra vt\ =\ \xi\hat\ra usvt\ =\ \xi\hat\ra u(s\la v)\tau(s\ra v,t)\, ((s\ra v)\cdot t)\ .
\]
It is sufficient to prove this with $u=e$, so we need to show
\begin{eqnarray}
((\xi\hat\ra s)\bar\ra v)\hat\ra t\ =\ (\xi\bar\ra (s\la v)\tau(s\ra v,t))\hat\ra ((s\ra v)\cdot t)\ .
\label{elp1}
\end{eqnarray}
By using the cross relation (\ref{ten2}) we see that, for $w\in G$ and $p\in M$,
\begin{eqnarray}
(\eta\bar\ra w)\hat\ra p\ =\ ((p^L\ra w^{-1}|\eta|^{-1})\bar\la \eta)\bar\ra (p^L\la w^{-1})^{-1}
\tau(p^L,p)\ .\label{elp3}
\end{eqnarray}
Using (\ref{elp3}) we calculate
\begin{eqnarray*}
((\xi\hat\ra s)\bar\ra v)\hat\ra t &=& 
(((s^L\ra|\xi|^{-1})\bar\la\xi)\bar\ra \tau(s^L,s)v)\hat\ra t  \cr
&=& ((t^L\ra v^{-1}\tau(s^L,s)^{-1}z^{-1})\bar\la
((s^L\ra|\xi|^{-1})\bar\la\xi))\ \bar\ra\  (t^L\la v^{-1}\tau(s^L,s)^{-1})^{-1}\tau(t^L,t)\ ,
\end{eqnarray*}
where $z=|(s^L\ra|\xi|^{-1})\bar\la\xi|$. Now, from (\ref{ten1}),
\[
(t^L\ra v^{-1}\tau(s^L,s)^{-1}z^{-1})\bar\la
((s^L\ra|\xi|^{-1})\bar\la\xi)\ =\ ((p'\cdot (s^L\ra|\xi|^{-1}))\bar\la\xi)\bar\ra
\tau(p'\ra((s^L\ra |\xi|^{-1})\la|\xi|),s^L)^{-1}\ ,
\]
where 
\begin{eqnarray*}
p' &=& t^L\ra v^{-1}\, \tau(s^L,s)^{-1}z^{-1}\,\tau(\langle (s^L\ra|\xi|^{-1})
\bar\la\xi\rangle ,s^L)\ \tau(s^L\ra
|\xi|^{-1},\langle \xi\rangle )^{-1}\cr &=&t^L\ra v^{-1}\tau(s^L,s)^{-1}
((s^L\ra|\xi|^{-1})\la|\xi|)^{-1}\ .
\end{eqnarray*}
From this we can calculate
\[
(p'\cdot(s^L\ra|\xi|^{-1}))\ra|\xi|\ =\ (t^L\ra v^{-1}\tau(s^L,s)^{-1})\cdot s^L\ ,
\]
so
\begin{eqnarray*}
((\xi\hat\ra s)\bar\ra v)\hat\ra t &=& ((((t^L\ra v^{-1}\tau(s^L,s)^{-1})\cdot s^L)
\ra|\xi|^{-1})\bar\la\xi)   \cr &&\quad
\bar\ra\tau(t^L\ra v^{-1}\tau(s^L,s)^{-1},s^L)^{-1}(t^L\la v^{-1}\tau(s^L,s)^{-1})^{-1}\tau(t^L,t)\ ,
\end{eqnarray*}
and set this equal to $((a\ra|\xi|^{-1})\bar\la\xi)\bar\ra y$. Now consider the right hand side of
(\ref{elp1}), using (\ref{elp3}):
\begin{eqnarray*}
 (\xi\bar\ra (s\la v)\tau(s\ra v,t))\hat\ra ((s\ra v)\cdot t) &=&
((((s\ra v)\cdot t)^L\ra\tau(s\ra v,t)^{-1}(s\la v)^{-1}|\xi|^{-1})\bar\la\xi) \cr
&& \bar\ra(((s\ra v)\cdot t)^L\la \tau(s\ra v,t)^{-1}(s\la v)^{-1})^{-1}
\tau(((s\ra v)\cdot t)^L,(s\ra v)\cdot t)\ ,
\end{eqnarray*}
which we set equal to $((b\ra|\xi|^{-1})\bar\la\xi)\bar\ra x$. 
It is our job to show that $a=b$ and $x=y$. We use the result,
derived from the equations $(c\cdot d)^L\,\tau(c,d)^{-1}cd=e$ and $c^{-1}=\tau(c^L,c)^{-1}c^L$,
\[
(c\cdot d)^L\ra\tau(c,d)^{-1}\ =\ (d^L\ra\tau(c^L,c)^{-1})\cdot c^L
\]
to show that
\[
b\ =\ ((t^L\ra\tau((s\ra v)^L,s\ra v)^{-1})\cdot (s\ra v)^L)\ \ra\  (s\la v)^{-1}\ .
\]
Now $((s\ra v)^L\cdot (s\ra v))\ra v^{-1}=e$, so $((s\ra v)^L\ra (s\la v)^{-1})\cdot s=e$, i.e.\
$(s\ra v)^L\ra (s\la v)^{-1}=s^L$.  On the other hand,
\begin{eqnarray*}
(s\ra v)^L\la (s\la v)^{-1} &=& (s\ra v)^L\la ((s\ra v)\la v^{-1})\ =\ \tau((s\ra v)^L,s\ra v)v^{-1}
\tau((s\ra v)^L\ra((s\ra v)\la v^{-1}),s)^{-1} \cr
&=& \tau((s\ra v)^L,s\ra v)v^{-1}\tau(s^L,s)^{-1}\ ,
\end{eqnarray*}
and we deduce that $a=b$. Next consider
\begin{eqnarray*}
b^{-1}x &=& (s\la v)\tau(s\ra v,t)\, ((s\ra v)\cdot t)^{L-1}\tau(((s\ra v)\cdot t)^L,(s\ra v)\cdot t)
\cr &=& (s\la v)\tau(s\ra v,t)\, ((s\ra v)\cdot t)\ =\  (s\la v)(s\ra v)t\ =\ svt\ , \cr
a^{-1}y &=& s^{L-1} (t^L\ra v^{-1}\tau(s^L,s)^{-1})^{-1}(t^L\la v^{-1}\tau(s^L,s)^{-1})^{-1}
\tau(t^L,t)\cr &=&  s^{L-1}\tau(s^L,s)\, v\, t^{L-1}\tau(t^L,t)\ =\ svt\ ,
\end{eqnarray*}
so we deduce that $x=y$.\endproof

\begin{propos} The $X$-action on tensor products in ${\cal D}$ is given by
\[
(\xi\tens\eta)\hat\ra x\ =\ (\xi\hat\ra (\|\eta\|\hat\la x))\tens \eta\hat\ra x\ .
\]
\end{propos}
\proof   We begin with
\begin{eqnarray*}
(\xi\tens\eta)\hat\ra u &=&(\xi\tens\eta)\bar\ra u\ =\ (\xi\bar\ra(\langle\eta\rangle\la u)
\tens \eta\bar\ra u\ ,
\end{eqnarray*}
so we consider $\|\eta\|\tilde\la u=\langle \eta\rangle u t^{-1}$ where $\|\eta\|\tilde\ra u=vt$, 
i.e.\ $t=\langle \eta\rangle \ra u$. Then we see that $\|\eta\|\tilde\la u=\langle \eta\rangle \la u$
as required. 
Continue with
\begin{eqnarray*}
(\xi\tens\eta)\hat\ra s\tau(s^L,s)^{-1} &=&   (\xi\tens\eta)\hat\ra s^{L-1} \ =\ 
 ((s^L\ra|\xi\tens\eta|^{-1})\bar\la(\xi\tens\eta) \cr
&=& ((s^L\ra |\eta|^{-1}|\xi|^{-1})\bar\la\xi)\bar\ra\tau(s^L\ra |\eta|^{-1},\langle \eta\rangle)
\tau(\langle (s^L\ra |\eta|^{-1})\bar\la\eta\rangle,s^L)^{-1}\cr &&\quad
\tens (s^L\ra |\eta|^{-1})\bar\la\eta\ .
\end{eqnarray*}
Now we can write
\[
(s^L\ra |\eta|^{-1}|\xi|^{-1})\bar\la\xi\ =\ \xi\hat\ra p\tau(p^L,p)^{-1}\ =\ \xi\hat\ra p^{L-1}\ ,
\]
where $p^L=s^L\ra|\eta|^{-1}$, and rewrite
\[
(\xi\tens\eta)\hat\ra s\tau(s^L,s)^{-1}\ =\ 
\xi\hat\ra p^{L-1}\tau(p^L,\langle \eta\rangle)
\tau(\langle p^L\bar\la\eta\rangle,s^L)^{-1}\tens (s^L\ra |\eta|^{-1})\bar\la\eta\ .
\]
We know that
\[
|p^L\bar\la\eta|^{-1}\langle p^L\bar\la\eta\rangle\ =\ 
(p^L\ra|\eta|)|\eta|^{-1}\langle \eta\rangle(p^L\ra|\eta|)^{-1}\ =\ 
s^L|\eta|^{-1}\langle \eta\rangle s^{L-1}\ =\ |\eta|^{-1}\langle \eta\rangle \hat\ra s^{L-1}\ ,
\]
 so
\begin{eqnarray*}
|\eta|^{-1}\langle \eta\rangle \hat\la s^{L-1} &=& \langle \eta\rangle s^{L-1}
\langle p^L\bar\la\eta\rangle^{-1}  \cr &=&
\langle \eta\rangle  (\langle p^L\bar\la\eta\rangle \cdot  s^{L})^{-1} \tau(
\langle p^L\bar\la\eta\rangle ,  s^{L})^{-1}\cr &=&
\langle \eta\rangle  (\langle p^L\bar\la\eta\rangle \cdot  (p^L\ra|\eta|))^{-1} \tau(
\langle p^L\bar\la\eta\rangle ,  s^{L})^{-1}\cr&=&
\langle \eta\rangle  ( p^L \cdot  \langle\eta\rangle)^{-1} \tau(
\langle p^L\bar\la\eta\rangle ,  s^{L})^{-1}\cr&=&
\langle \eta\rangle \langle\eta\rangle^{-1}  p^{L-1} \tau(  p^L, \langle\eta\rangle) \tau(
\langle p^L\bar\la\eta\rangle ,  s^{L})^{-1}\ ,
\end{eqnarray*}
as required\endproof

\begin{propos} The braiding $\Psi$ is given in terms of the $X$-action by
\[
\Psi(\xi\tens\eta)\ =\ \eta\hat\ra(\langle\xi\rangle\ra|\eta|)^{-1}\tens\xi\hat\ra|\eta|\quad,\quad
\Psi^{-1}(\xi'\tens\eta')\ =\ \eta'\hat\ra| \xi'\hat\ra\langle\eta'\rangle
 |^{-1} \tens \xi'\hat\ra\langle\eta'\rangle\ .
\]
\end{propos}
\proof
   By definition of $\hat\ra s$,
\[
\eta\hat\ra s\tau(s^L,s)^{-1}\ =\ \eta\hat\ra s^{L-1}\ =\ (s^L\ra|\eta|^{-1})\bar\la\eta\ ,
\]
so we deduce that $t\bar\la \eta=\eta\hat\ra(t\ra|\eta|)^{-1}$. Now
\[
\Psi(\xi\tens\eta)\ =\ \langle \xi\rangle \bar\la\eta\tens\xi\bar\ra|\eta|\ =\ 
\eta\hat\ra (\langle \xi\rangle \ra|\eta|)^{-1}\tens \xi\hat\ra|\eta|\ .\square
\]

\section{A bialgebra in the braided category}
Take a group $X$ with subgroup $G$, and
a set of left coset representatives $M$ which contains $e$. 
We assume that $(M,\cdot)$ has the left division property, 
i.e.\  for all $t,s\in M$ there is a unique solution $p\in M$ to the equation $s\cdot p=t$.

Introduce a vector space $D$ with basis $\delta_y\tens x$
for $y\in Y$ and $x\in X$. Then we define
\[
\xi \hat\ra(\delta_y\tens x)\ =\ \delta_{y,\|\xi\|}\ \xi\hat\ra x\ .
\]
We see that $D$ is an object of ${\cal D}$ with grade $y\circ\|\delta_y\tens x\|=y\tilde\ra x$
and action
\[
(\delta_y\tens x)\hat\ra z\ =\ \delta_{y\tilde\ra(a\tilde\la z)}\tens (a\tilde\la z)^{-1}xz\ ,
\]
where $a=\|\delta_y\tens x\|$. Then the associative
 multiplication $\mu$ on $D$ consistent with the action
is 
\[
(\delta_y\tens x)(\delta_w\tens z)\ =\ \delta_{w,y\tilde\ra x}\ \delta_{y\tilde\ra 
\tilde\tau(a,b)}\tens \tilde\tau(a,b)^{-1}xz\ ,
\]
where $b=\|\delta_w\tens z\|$, and is a morphism in ${\cal D}$. 

This much we have done before
in ${\cal C}$. The additional ingredient we have in ${\cal D}$ is the braiding. 
We can use the braiding to define a coproduct
for $D$ which in turn gives the tensor product structure in ${\cal D}$. 

\begin{propos} The coproduct in $D$ consistent with the tensor product structure in ${\cal D}$
is
\begin{eqnarray}
\Delta (\delta_y\tens x) &=&
\sum_{z,w\in Y: z w=y} 
\delta_{w\tilde \ra   |h||h_{(2)}|^{-1}|h_{(1)}|^{-1}} \tens 
|h_{(1)}|\,|h_{(2)}|\, |h|^{-1}
 x\,  \langle h_{(2)}\rangle^{-1} \cr &&\qquad \tens 
\delta_{z\tilde\ra |h||h_{(2)}|^{-1}}\tens
|h_{(2)}||h|^{-1} x\ ,
\end{eqnarray}
where $h=\delta_y\tens x$, $
\|h_{(2)}\|= |h|^{-1}z^{-1}|h|\, x^{-1}z\, x
$, and  $\|h_{(1)}\|=|h_{(2)}|\, \|h\|\, \<h_{(2)}\>^{-1}$.
\end{propos}
\proof Begin with
\begin{eqnarray}
(\xi\tens\eta)\hat\ra(\delta_y\tens x)\ =\ \delta_{y,\|\xi\|\circ\|\eta\|}\ (\xi\tens\eta)\hat\ra x
\ =\ \delta_{y,\|\xi\|\circ\|\eta\|}\ \xi\hat\ra(\|\eta\|\tilde\la x)\tens \eta\hat\ra x\ .
\label{ttt1}
\end{eqnarray}
For $\Delta h=h_{(1)}\tens h_{(2)}$ this should be the same as (using 
$\|h\|=\|h_{(1)}\|\circ \|h_{(2)}\|$)
\begin{eqnarray}
(\xi\tens\eta)\tens(h_{(1)}\tens h_{(2)}) &\stackrel{\Phi}\longmapsto &
\xi\hat\ra\tilde\tau(\|\eta\|,\|h\|)\tens(\eta\tens(h_{(1)}\tens h_{(2)})) \cr
 &\stackrel{I\tens\Phi^{-1}}\longmapsto &
\xi\hat\ra\tilde\tau(\|\eta\|,\|h\|)\tens((\eta\hat\ra\tilde\tau(\|h_{(1)}\|,\|h_{(2)}\|)^{-1}
\tens h_{(1)})\tens h_{(2)})\cr
 &\stackrel{I\tens(\Psi\tens I)}\longmapsto &
\xi\hat\ra\tilde\tau(\|\eta\|,\|h\|)\tens((h'_{(1)}
\tens \eta')\tens h_{(2)})
\cr
 &\stackrel{I\tens\Phi}\longmapsto &
\xi\hat\ra\tilde\tau(\|\eta\|,\|h\|)\tens(h'_{(1)}\hat\ra\tilde\tau(\|\eta'\|,\|h_{(2)}\|)
\tens (\eta'\tens h_{(2)}))\cr
 &\stackrel{\Phi^{-1}}\longmapsto &
(\xi\hat\ra\tilde\tau(\|\eta\|,\|h\|)n\tens h'_{(1)}\hat\ra\tilde\tau(\|\eta'\|,\|h_{(2)}\|))
\tens (\eta'\tens h_{(2)})\cr
 &\stackrel{\hat\ra\tens\hat\ra}\longmapsto &
(\xi\hat\ra\tilde\tau(\|\eta\|,\|h\|)n)\hat\ra( h'_{(1)}\hat\ra\tilde\tau(\|\eta'\|,\|h_{(2)}\|))
\tens \eta'\hat\ra h_{(2)}\ ,
\label{ttt2}\end{eqnarray}
where 
\begin{eqnarray*}
n &=& \tilde\tau(\| h'_{(1)}\hat\ra\tilde\tau(\|\eta'\|,\|h_{(2)}\|)  \|, \|  \eta' \|\circ\|  h_{(2)} \|)^{-1}
\cr &=&
 \tilde\tau(\| h'_{(1)}\|\tilde\ra\tilde\tau(\|\eta'\|,\|h_{(2)}\|)  , \|  \eta' \|\circ\|  h_{(2)} \|)^{-1}\ ,
   \cr h'_{(1)}\tens \eta' &=& \Psi(\eta\hat\ra\tilde\tau(\|h_{(1)}\|,\|h_{(2)}\|)^{-1}
\tens h_{(1)}) \cr &=&
h_{(1)} \hat\ra ( \langle\eta\hat\ra\tilde\tau(\|h_{(1)}\|,\|h_{(2)}\|)^{-1}\rangle \ra |  h_{(1)} |   )^{-1}
\tens \eta\hat\ra\tilde\tau(\|h_{(1)}\|,\|h_{(2)}\|)^{-1} |h_{(1)}|
 \ .
\end{eqnarray*}
As $\tilde\tau$ takes values in $G$, we see $h'_{(1)}=h_{(1)}\hat\ra \langle \eta'\rangle^{-1}$. 
If we set
\begin{eqnarray}
h_{(1)}\hat\ra \langle \eta'\rangle^{-1}\tilde\tau(\|\eta'\|,\|h_{(2)}\|) &=&   
\delta_{\|\xi\|\tilde \ra\tilde\tau(\|\eta\|,\|h\|)n} \tens
n^{-1}\tilde\tau(\|\eta\|,\|h\|)^{-1}(\|\eta\|\tilde\la x)\ =\ h'\ ,\cr
h_{(2)} &=&\delta_{\|\eta\|\tilde\ra \tilde\tau(\|h_{(1)}\|,\|h_{(2)}\|)^{-1} |h_{(1)}|}\tens
|h_{(1)}|^{-1} \tilde\tau(\|h_{(1)}\|,\|h_{(2)}\|) x\ ,\label{ttt3}
\end{eqnarray}
then (\ref{ttt1}) and (\ref{ttt2}) agree if $y=\|\xi\|\circ\|\eta\|$. 
Now define $vt$ ($v\in G$ and $t\in M$) by
the factorisation $\|\eta\|\tilde\ra x=\|\eta'\|\circ \|h_{(2)} \|=vt$,
and from this $t=\langle\eta'\rangle\cdot \langle h_{(2)}\rangle$ and
$\|\eta\|\tilde\la x=\langle\eta\rangle x t^{-1}$. Also we have
$\langle \eta'\rangle^{-1}\tau(\langle \eta'\rangle,\langle h_{(2)}\rangle)
=\langle h_{(2)}\rangle (\langle \eta'\rangle\cdot\langle h_{(2)}\rangle)^{-1}=
\langle h_{(2)}\rangle t^{-1}$ and $n=\tau(\langle h'\rangle,t)^{-1}$. Now
$\|h'\|\tilde\la t=\langle h'\rangle tp^{-1}$
where $\|h'\|\tilde\ra t=\|h'\hat\ra t\|=up$ ($u\in G$ and $p\in M$). Then
\begin{eqnarray*}
\tilde\tau(\|\eta\|,\|h\|)n(\|h'\|\tilde\la t) &=
&\tau(\langle\eta\rangle,\langle h\rangle)\,n\,\langle h'\rangle\, t\,p^{-1}
\cr &=& \tau(\langle\eta\rangle,\langle h\rangle)(\langle h'\rangle\cdot t)p^{-1}
\cr &=& \tau(\langle\eta\rangle,\langle h\rangle)(\langle \eta\rangle\cdot \langle h\rangle)p^{-1}\ =\ 
\langle \eta\rangle\, \langle h\rangle\, p^{-1}\ ,
\end{eqnarray*}
(using $\|h'\|\circ t=\|\eta\|\circ \|h\|$ from (\ref{ttt2})), so from  (\ref{ttt3}),
\begin{eqnarray}
h_{(1)}\hat\ra \langle h_{(2)}\rangle &=& h'\hat\ra t\ =\ 
\delta_{\|\xi\|\tilde \ra\tilde\tau(\|\eta\|,\|h\|)n(\|h'\|\tilde\la t)} \tens (\|h'\|\tilde\la t)^{-1}
n^{-1}\tilde\tau(\|\eta\|,\|h\|)^{-1}(\|\eta\|\tilde\la x)t \cr &=&
\delta_{\|\xi\|\tilde \ra\langle \eta\rangle\,  \langle h\rangle\, p^{-1}} \tens p\, \langle h\rangle^{-1}
 x\ .\label{ttt9}
\end{eqnarray}
Now we calculate, using $\|h_{(1)}\|\tilde\ra \langle h_{(2)}\rangle = \|h'\|\tilde\ra t=up$,
\[
\| h_{(1)} \hat\ra \langle h_{(2)}\rangle\| \tilde\la \langle h_{(2)}\rangle^{-1}\ =\ 
(\| h_{(1)}\| \tilde\la \langle h_{(2)}\rangle)^{-1}\ =\ (\langle h_{(1)}\rangle\langle h_{(2)}\rangle
 p^{-1})^{-1}\ ,
\]
so if we apply $\hat\ra \langle h_{(2)}\rangle^{-1}$ to (\ref{ttt9}) we get
\begin{eqnarray}
h_{(1)} &=& 
\delta_{\|\xi\|\tilde \ra\langle \eta\rangle\,  \langle h\rangle
\, (\langle h_{(1)}\rangle\langle h_{(2)}\rangle)^{-1}} \tens 
\langle h_{(1)}\rangle\langle h_{(2)}\rangle\, \langle h\rangle^{-1}
 x\,  \langle h_{(2)}\rangle^{-1}\cr
&=& 
\delta_{\|\xi\|\tilde \ra\langle\eta\rangle\,  \tau(\langle h_{(1)}\rangle,\langle h_{(2)}\rangle)^{-1}} \tens 
\tau(\langle h_{(1)}\rangle,\langle h_{(2)}\rangle)\, 
 x\,  \langle h_{(2)}\rangle^{-1}
\ .\label{ttt10}
\end{eqnarray}
 From 
$\|h\|=\|h_{(1)}\|\circ \|h_{(2)}\|$ we see that $\tau(\langle h_{(1)}\rangle,
\langle h_{(2)}\rangle)^{-1}|h_{(1)}|=|h||h_{(2)}|^{-1}$, so we can rewrite (\ref{ttt3}) to give
\begin{eqnarray}
h_{(2)} &=&\delta_{\|\eta\|\tilde\ra |h||h_{(2)}|^{-1}}\tens
|h_{(2)}||h|^{-1} x\ .\label{ttt5}
\end{eqnarray}
Now we use the definition of $\|h_{(2)}\|$ on this formula to get
\[
(\|\eta\|\tilde\ra |h||h_{(2)}|^{-1})\circ \|h_{(2)}\|\ =\ |h|^{-1}\|\eta\|\, |h|\,|h_{(2)}|^{-1}\langle h_{(2)}\rangle
\ =\ \|\eta\|\tilde\ra x\ ,
\]
which we rearrange as
$
\|h_{(2)}\|= |h|^{-1}\|\eta\|^{-1}|h|\, x^{-1}\|\eta\|\, x
$. If we set $\|\eta\|=z$ and $\|\xi\|=z'$, we 
have the constraint $\|\xi\|\circ\|\eta\|=y=|\eta|^{-1}z'\<\eta\>=z(z'\tilde\ra \<\eta\>)$, and if we set
$w=z'\tilde\ra \<\eta\>$ then
\begin{eqnarray}
\Delta (\delta_y\tens x) &=&
\sum_{z,w\in Y: z w=y} 
\delta_{w\tilde \ra   \tau(\langle h_{(1)}\rangle,\langle h_{(2)}\rangle)^{-1}} \tens 
\tau(\langle h_{(1)}\rangle,\langle h_{(2)}\rangle)
 x\,  \langle h_{(2)}\rangle^{-1} \cr &&\quad \tens 
\delta_{z\tilde\ra |h||h_{(2)}|^{-1}}\tens
|h_{(2)}||h|^{-1} x\ ,
\end{eqnarray}
where $
\|h_{(2)}\|= |h|^{-1}z^{-1}|h|\, x^{-1}z\, x
$. We find $\|h_{(1)}\|$ by solving the equation $\|h_{(1)}\|\circ \|h_{(2)}\|=\|h\|$
to get $\|h_{(1)}\|=|h_{(2)}|\, \|h\|\, \<h_{(2)}\>^{-1}$. Finally we can substitute
$\tau(\langle h_{(1)}\rangle,\langle h_{(2)}\rangle)=|h_{(1)}|\,|h_{(2)}|\, |h|^{-1}$. 

\endproof

\begin{propos} The map $\epsilon:D\to k$ defined by 
$\epsilon(\delta_y\tens x)=\delta_{y,e}$ is a counit for the coproduct,
and $\Delta(I)=I\tens I$.
\end{propos}
\proof  Let $h=\delta_y\tens x$. For $(\epsilon\tens I)\Delta$ begin with
\begin{eqnarray}
(\epsilon\tens I)\Delta (\delta_y\tens x) &=&
\sum_{z,w\in Y: z w=y}  \delta_{w,e}\ 
\delta_{z\tilde\ra |h||h_{(2)}|^{-1}}\tens
|h_{(2)}||h|^{-1} x \cr
&=& \delta_{y\tilde\ra |h||h_{(2)}|^{-1}}\tens
|h_{(2)}||h|^{-1} x\ ,
\end{eqnarray}
where $\|h_{(2)}\|= |h|^{-1}y^{-1}|h|\, x^{-1}y\, x$. But by definition
$y\circ\|h\|=|h|^{-1}y\<h\>=x^{-1}yx$, so we deduce that $\|h_{(2)}\|=\|h\|$
and $\|h_{(1)}\|=e$.

Now for $(I\tens\epsilon)\Delta$ begin with
\begin{eqnarray}
(I\tens\epsilon)\Delta (\delta_y\tens x)\ =\ 
\sum_{z,w\in Y: z w=y}   \delta_{z,e}\ 
\delta_{w\tilde \ra   |h||h_{(2)}|^{-1}|h_{(1)}|^{-1}} \tens 
|h_{(1)}|\,|h_{(2)}|\, |h|^{-1}
 x\,  \langle h_{(2)}\rangle^{-1} \ ,
\end{eqnarray}
where  $\|h_{(2)}\|= |h|^{-1}z^{-1}|h|\, x^{-1}z\, x=e$
and $\|h_{(1)}\|=|h_{(2)}|\, \|h\|\, \<h_{(2)}\>^{-1}=\|h\|$. 

The proof that $\Delta(I)=I\tens I$ is easy once you notice that
for every $h=\delta_y\tens x$ term in $I$, we have $x=\|h\|=e$. \endproof

As the formula for the coproduct is not very nice, we shall use standard diagramatic arguments
to show that $(D,\mu,\Delta)$ is a bialgebra \cite{Ma:book}. The pentagon identity means that
we do not have to keep track of every re-bracketing done in the course of following a diagram.
We just assume that there is a fixed bracketing at the beginning and at the end,
and apply the associator as required in between. 
In Fig 1 we give
(in order) the symbols we shall use for
the braiding $\Psi$, the action $\hat\ra:V\tens D\to V$, the counit, unit, product and coproduct:
$$
\epsfig{file=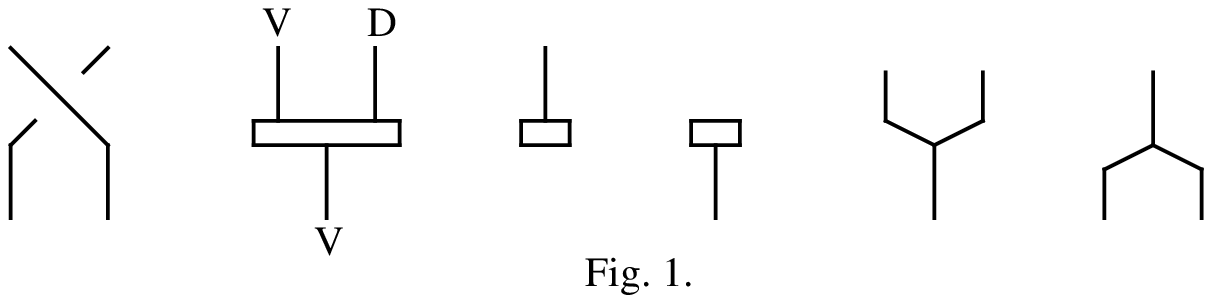}
$$
%\psboxscaled{1000}{fig1.ps}%{Fig\ 2}
The similarity between the symbols for the action and the counit is not coincidental.
The counit is the action of $D$ on $k$, which is traditionally represented by an invisible line. 
In Fig 2 we give the {\bf definition} of the product $\mu$ and coproduct $\Delta$ on $D$. 
$$
\epsfig{file=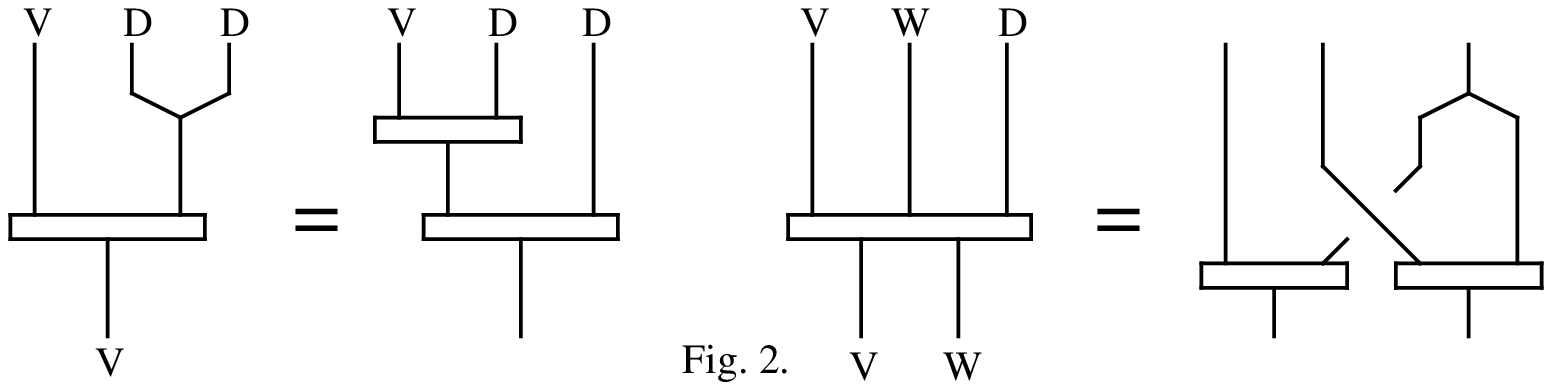}
$$
%\psboxscaled{1000}{fig2.ps}%{Fig\ 2}
Now the proof that $\Delta$ is multiplicative is given as: (Fig 3)
$$
\epsfig{file=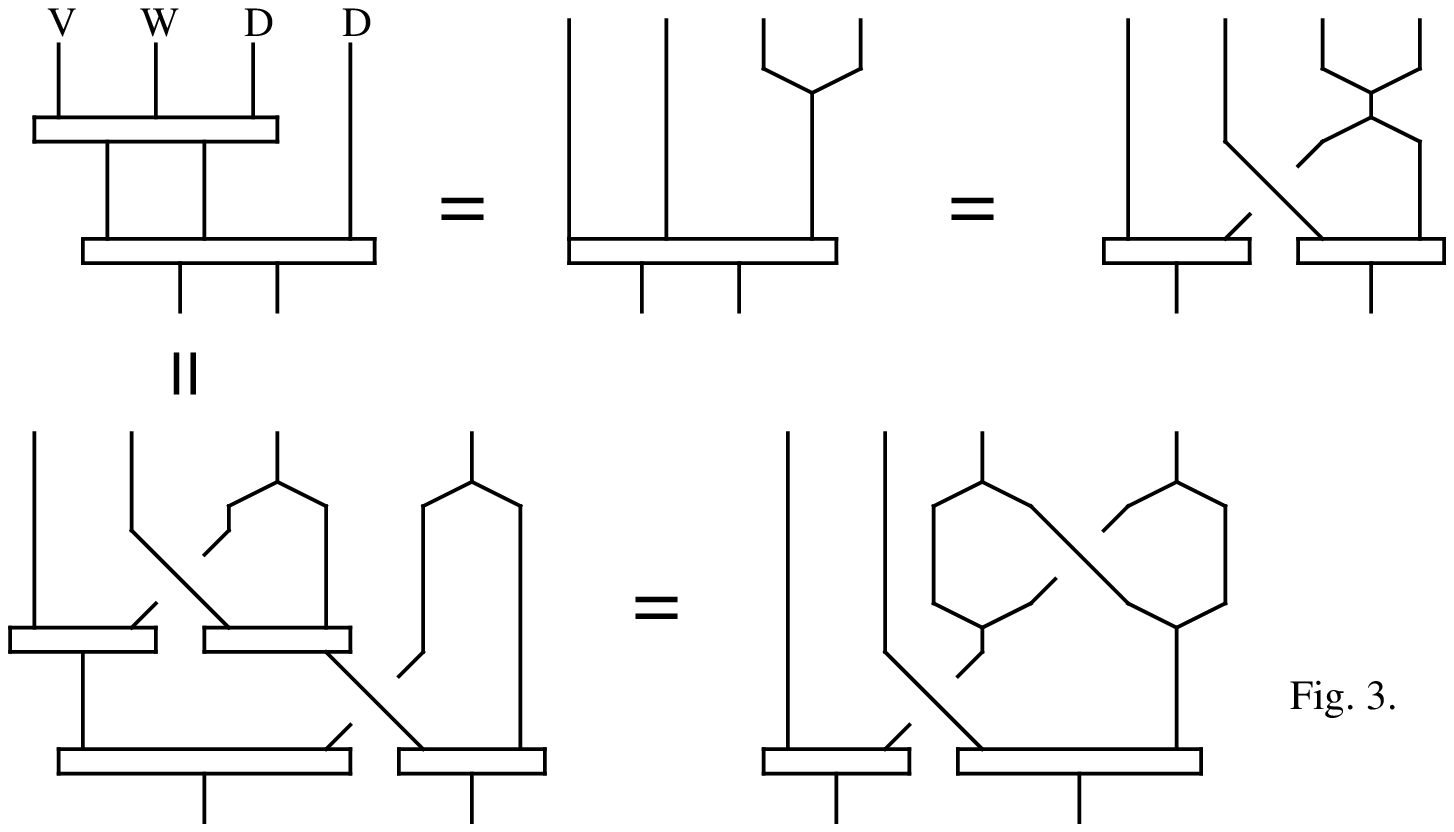}
$$
%\psboxscaled{1000}{fig3.ps}%{Fig\ 2}
To show that $\Delta$ is coassociative we must first show
$(\Phi((\xi\tens\eta)\tens\kappa))\hat\ra h=\Phi(((\xi\tens\eta)\tens\kappa)\hat\ra h)$, 
which is easy enough to check from the definitions. This then means that the following two
ways of splitting up the calculation of the action on a triple tensor product are the same: (Fig 4)
$$
\epsfig{file=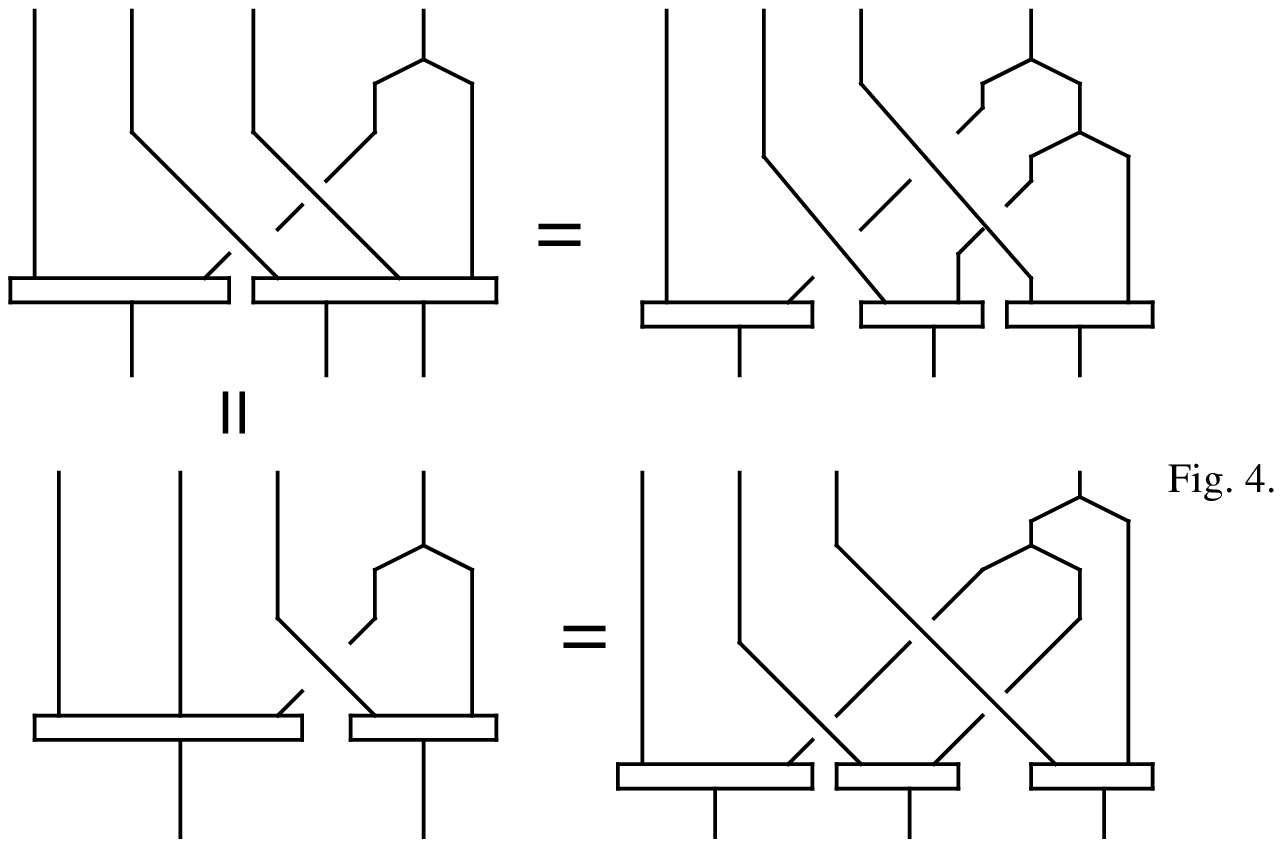}
$$
%\psboxscaled{1000}{fig4.ps}%{Fig\ 2}

\section{A rigid braided tensor category}
We assume the same conditions on $(M,\cdot)$ as the last section.
Note that $(Y\circ)$ then has right inverses.
The definitions of dual, and the corresponding evaluation and coevaluation maps,
considered previously for ${\cal C}$, can also be used in ${\cal D}$. Fig 5(a) and 5(b) show the
diagrams we shall use for evaluation and coevaluation. Recall that the morphisms in 
${\cal D}$ are required to preserve the actions and gradings. This means that
if $T:V\to W$ is a morphism, then we have the picture in Fig 5(c). 
$$
\epsfig{file=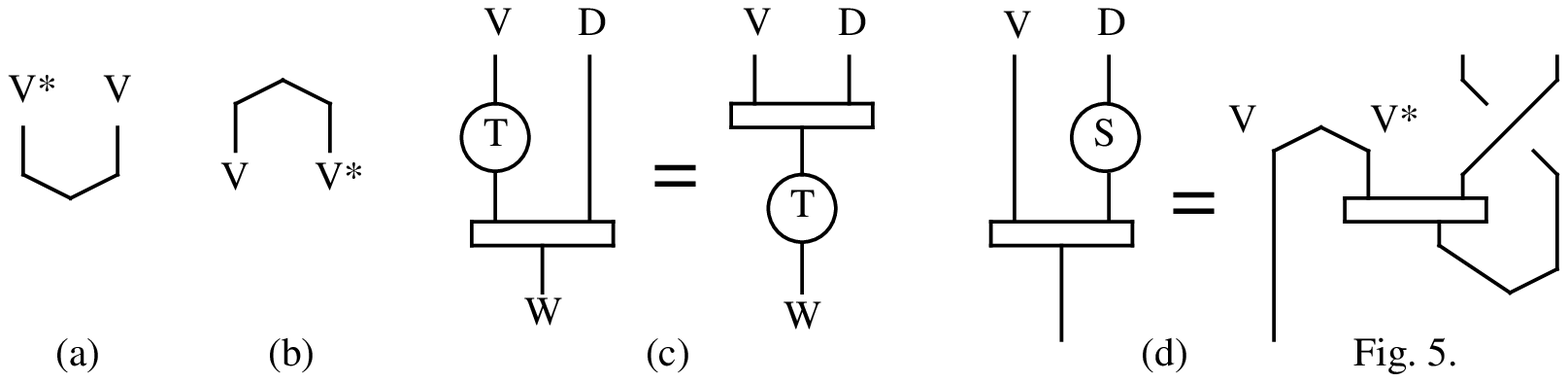}
$$
%\psboxscaled{1000}{fig5.ps}%{Fig\ 2}
We would like to show that $D$ is a braided Hopf algebra in the category ${\cal D}$,
and in Fig 5(d) we
give the {\bf definition} of antipode $S:D\to D$. Note that this is not the same picture 
as that in \cite{Ma:book}. This is because we are using right actions instead of
left actions. We cannot simply reflect the picture in \cite{Ma:book} either, as
the evaulation and coevaluation morphisms have a definite handedness. 
In the next proposition we find what the formula
for the antipode actually is, and then we go on to check that $S$ satisfies the 
required condition in the axioms of a braided Hopf algebra.

\begin{propos} Let $h=\delta_y\tens x\in D$. Then
$
S(h)\ =\ \delta_{y^{-1}|h|\langle h\rangle^{-1}}\tens \langle h\rangle
\,x^{-1}\, |h|
$.
\end{propos}
\proof
Suppose that $S(h)=\delta_{y'}\tens x'$.  Then 
for $\xi\in V$ we have $\xi\hat\ra S(h)=\delta_{y',\|\xi\|}\xi\hat\ra x'$, which by definition
 is equal to the composition
\begin{eqnarray}
\xi\tens h &\stackrel{{\rm coeval}_V\tens \Psi^{-1}}
\longmapsto&\sum (\eta\hat\ra\tau(\langle\eta\rangle^L,
\langle\eta\rangle)^{-1} \tens\hat\eta)\tens (h\hat\ra|\xi\hat\ra\langle h\rangle|^{-1}
\tens\xi\hat\ra \langle h\rangle) \cr
&\stackrel{ \Phi}\longmapsto&\sum \eta\hat\ra\tau(\langle\eta\rangle^L,
\langle\eta\rangle)^{-1}\tau(\langle\eta\rangle^L,\langle\xi\rangle\cdot\langle h\rangle)
 \tens(\hat\eta\tens (h\hat\ra|\xi\hat\ra\langle h\rangle|^{-1}\tens\xi\hat\ra \langle h\rangle)) \cr
&\stackrel{ I\tens\Phi^{-1}}\longmapsto&\sum \eta'
 \tens((\alpha 
\tens h\hat\ra|\xi\hat\ra\langle h\rangle|^{-1})\tens\xi\hat\ra \langle h\rangle) \cr
&\stackrel{ I\tens(\hat\ra\tens I)}\longmapsto&\sum \eta'
 \tens(\alpha 
\hat\ra (h\hat\ra|\xi\hat\ra\langle h\rangle|^{-1})\tens\xi\hat\ra \langle h\rangle) \cr
&\stackrel{ I\tens{\rm eval}}\longmapsto& \sum\eta'
 \ (\alpha 
\hat\ra (h\hat\ra|\xi\hat\ra\langle h\rangle|^{-1}))(\xi\hat\ra \langle h\rangle)
\ ,\label{1066}
\end{eqnarray}
where $\eta'=\eta\hat\ra\tau(\langle\eta\rangle^L,
\langle\eta\rangle)^{-1}\tau(\langle\eta\rangle^L,\langle\xi\rangle\cdot\langle h\rangle)$
and $\alpha=\hat\eta\hat\ra\tau(\langle h\hat\ra|\xi\hat\ra\langle h\rangle|^{-1}
\rangle,\langle\xi\hat\ra \langle h\rangle\rangle)^{-1}$. 
We set $u=|\xi\hat\ra\langle h\rangle|$, $s=\langle h\rangle \ra u^{-1}$
and $t=\langle \xi\hat\ra\langle h\rangle \rangle$. We have
\[
h\hat\ra u^{-1}\ =\ \delta_{y\tilde\ra(\|h\|\tilde\la u^{-1})}\tens 
(\|h\|\tilde\la u^{-1})^{-1}xu^{-1}\ .
\]
Now $\|h\|\tilde\ra u^{-1}=vs$, so $\|h\|\tilde\la u^{-1}=\langle h\rangle u^{-1}s^{-1}$, and
\[
(\hat\eta\hat\ra\tau(s,t)^{-1})\hat\ra(h\hat\ra u^{-1})\ =\ 
\delta_{\| \eta\|^L\tilde\ra \tau(s,t)^{-1},y\tilde\ra 
\langle h\rangle u^{-1}s^{-1}}\ \hat\eta\hat\ra 
\tau(s,t)^{-1} su\langle h\rangle^{-1}xu^{-1}\ .
\]
Now we recall that $\tau(s,t)^{-1} s=(s\cdot t)t^{-1}$, where $s\cdot t=
\langle \xi\rangle\cdot \langle h\rangle$ as $\Psi^{-1}$ preserves grades. Then
\[
(\hat\eta\hat\ra\tau(s,t)^{-1})\hat\ra(h\hat\ra u^{-1})\ =\ 
\delta_{y,\| \eta\|^L\tilde\ra (\langle \xi\rangle\cdot \langle h\rangle)t^{-1}u\langle h\rangle^{-1}
}\ \hat\eta\hat\ra 
(\langle \xi\rangle\cdot \langle h\rangle)t^{-1}u\langle h\rangle^{-1}xu^{-1}\ .
\]
Now $u^{-1}t=\|\xi\hat\ra\langle h\rangle\|=\|\xi\|\tilde\ra \langle h\rangle$, so
\begin{eqnarray}
(\hat\eta\hat\ra\tau(s,t)^{-1})\hat\ra(h\hat\ra u^{-1})\ =\ 
\delta_{y,\| \eta\|^{LD}\tilde\ra \tau(\langle \xi\rangle, \langle h\rangle)^{-1}|\xi|
}\ \hat\eta\hat\ra 
\tau(\langle \xi\rangle, \langle h\rangle)^{-1}|\xi|xu^{-1}\ .\label{trub}
\end{eqnarray}
Set $w=\tau(\langle \xi\rangle, \langle h\rangle)^{-1}|\xi|xu^{-1}$ and
$z=\|\xi\|\tilde\ra\langle h\rangle$. If $(\hat\eta\hat\ra w)(\xi\hat\ra \langle h\rangle)$
is not zero, then $\|\hat\eta\|\tilde\ra w=z^{LD}$. Then $z^{LD}\tilde\ra w^{-1}=\|\hat\eta\|
=\|\eta\|^{LD}$, so $(z^{LD}\tilde\ra w^{-1})^{RD}=\|\eta\|$. If (\ref{trub}) is not zero, then
\begin{eqnarray}
y &=& z^{LD}\tilde\ra w^{-1}\,  \tau(\langle \xi\rangle, \langle h\rangle)^{-1}\,|\xi|
\cr
 &=&(ut^{-1})\tilde\ra u x^{-1}\ =\ (t^{-1}u)\tilde\ra  x^{-1}\ ,
\cr
y^{-1} &=& u^{-1}t\tilde\ra  x^{-1}\ =\ \|\xi\|\tilde\ra \<h\>\,  x^{-1}\ ,\label{uup}
\end{eqnarray}
so we deduce that $y'=\|\xi\|=y^{-1}\tilde\ra x\,\langle h\rangle^{-1}$. 
By definition $y\circ\|h\|=|h|^{-1}y\<h\>=x^{-1}yx$, so we get
 $y'=\|\xi\|=\langle h\rangle x^{-1}y^{-1}x\,\langle h\rangle^{-1}=y^{-1}|h|\langle h\rangle^{-1}$. 
Now, from (\ref{dact}),
\begin{eqnarray}
(\hat\eta\hat\ra w)(\xi\hat\ra \langle h\rangle)
&=&\hat\eta(\xi\hat\ra\langle h\rangle
\,\tilde\tau(z^{LD},z)^{-1}(z^{LD}\tilde\la w^{-1})\,\tilde\tau(z^{LD}\tilde\ra w^{-1},
(z^{LD}\tilde\ra w^{-1})^{RD}))
\cr&=&\hat\eta(\xi\hat\ra\langle h\rangle
\,\tau(t^{L},t)^{-1}(z^{LD}\tilde\la w^{-1})\,\tau(\<\eta\>^L,\<\eta\>)) 
\cr&=&\hat\eta(\xi\hat\ra\langle h\rangle
\,\tau(t^{L},t)^{-1}\, t^L\,w^{-1}\<\eta\>^{L-1}\,\tau(\<\eta\>^L,\<\eta\>)) 
\cr&=&\hat\eta(\xi\hat\ra\langle h\rangle
\,t^{-1}\, w^{-1}\,\<\eta\>) 
\cr&=&\hat\eta(\xi\hat\ra\langle h\rangle
\,t^{-1}\, u\,x^{-1}\,|\xi|^{-1}\,\tau(\<\xi\>,\<h\>)\,\<\eta\>) 
\cr&=&\hat\eta(\xi\hat\ra\langle h\rangle
\,x^{-1}\, y\,|\xi|^{-1}\,\tau(\<\xi\>,\<h\>)\,\<\eta\>) \ .
\end{eqnarray}
We choose a basis element for the summation to be
$\eta=\xi\hat\ra\langle h\rangle
\,x^{-1}\, y\,|\xi|^{-1}\,\tau(\<\xi\>,\<h\>)\,\<\eta\>$, and then
\begin{eqnarray}
\xi\hat\ra x' &=& \eta'\ =\ \xi\hat\ra\langle h\rangle
\,x^{-1}\, y\,|\xi|^{-1}\,\tau(\<\xi\>,\<h\>)\,\<\eta\>\, 
\tau(\langle\eta\rangle^L,
\langle\eta\rangle)^{-1}\tau(\langle\eta\rangle^L,\langle\xi\rangle\cdot\langle h\rangle)
\cr
&=&  \xi\hat\ra\langle h\rangle
\,x^{-1}\, y\,|\xi|^{-1}\,\tau(\<\xi\>,\<h\>)\,
\langle\eta\rangle^{L-1}
\,\tau(\langle\eta\rangle^L,\langle\xi\rangle\cdot\langle h\rangle) \cr
&=&  \xi\hat\ra\langle h\rangle
\,x^{-1}\, y\,|\xi|^{-1}\,\tau(\<\xi\>,\<h\>)\,
(\langle\xi\rangle\cdot\langle h\rangle)
\,(\langle\eta\rangle^L\cdot(\langle\xi\rangle\cdot\langle h\rangle) )^{-1}\cr
&=&  \xi\hat\ra\langle h\rangle
\,x^{-1}\, y\,|\xi|^{-1}\,\<\xi\>\,\<h\>
\,(((\langle\eta\rangle^L \ra\tau(\langle\xi\rangle,\langle h\rangle)^{-1})
\cdot \langle\xi\rangle)\cdot\langle h\rangle )^{-1}\cr
&=&  \xi\hat\ra\langle h\rangle
\,x^{-1}\, |h|
\,(((\langle\eta\rangle^L \ra\tau(\langle\xi\rangle,\langle h\rangle)^{-1})
\cdot \langle\xi\rangle)\cdot\langle h\rangle )^{-1}
\end{eqnarray}
Finally, from the top line of (\ref{uup}), we get $\|\eta\|^{LD}
\tilde\ra  \tau(\langle \xi\rangle, \langle h\rangle)^{-1}=y \tilde\ra |\xi|^{-1}$. Then 
if we set $c=(\langle\eta\rangle^L \ra\tau(\langle\xi\rangle,\langle h\rangle)^{-1})
\cdot \langle\xi\rangle$, we see that $vc=(y \tilde\ra |\xi|^{-1})\circ y'$ (for some $v\in G$). 
But
\[
(y \tilde\ra |\xi|^{-1})\circ y'\ =\ (|\xi|\,y\,|\xi|^{-1})\circ \|\xi\|\ =\ 
|\xi|^{-1}\, |\xi|\,y\,|\xi|^{-1}\,\<\xi\>\ =\ y\,\|\xi\|\ =\ |h|\,\<h\>^{-1}\ ,
\]
so $c=\<h\>^L$. We conclude that $x'=\langle h\rangle
\,x^{-1}\, |h|$. \endproof

\begin{propos} The morphisms $\mu(I\tens S)\Delta:D\to D$ and $\mu(S\tens I)\Delta:D\to D$
are both equal to $1.\epsilon:D\to D$. 
\end{propos}

\newpage
\proof This part of the definition of a braided Hopf algebra can be checked by diagrams.
First for $\mu(I\tens S)\Delta:D\to D$ we have Fig.\ 6.
$$
\epsfig{file=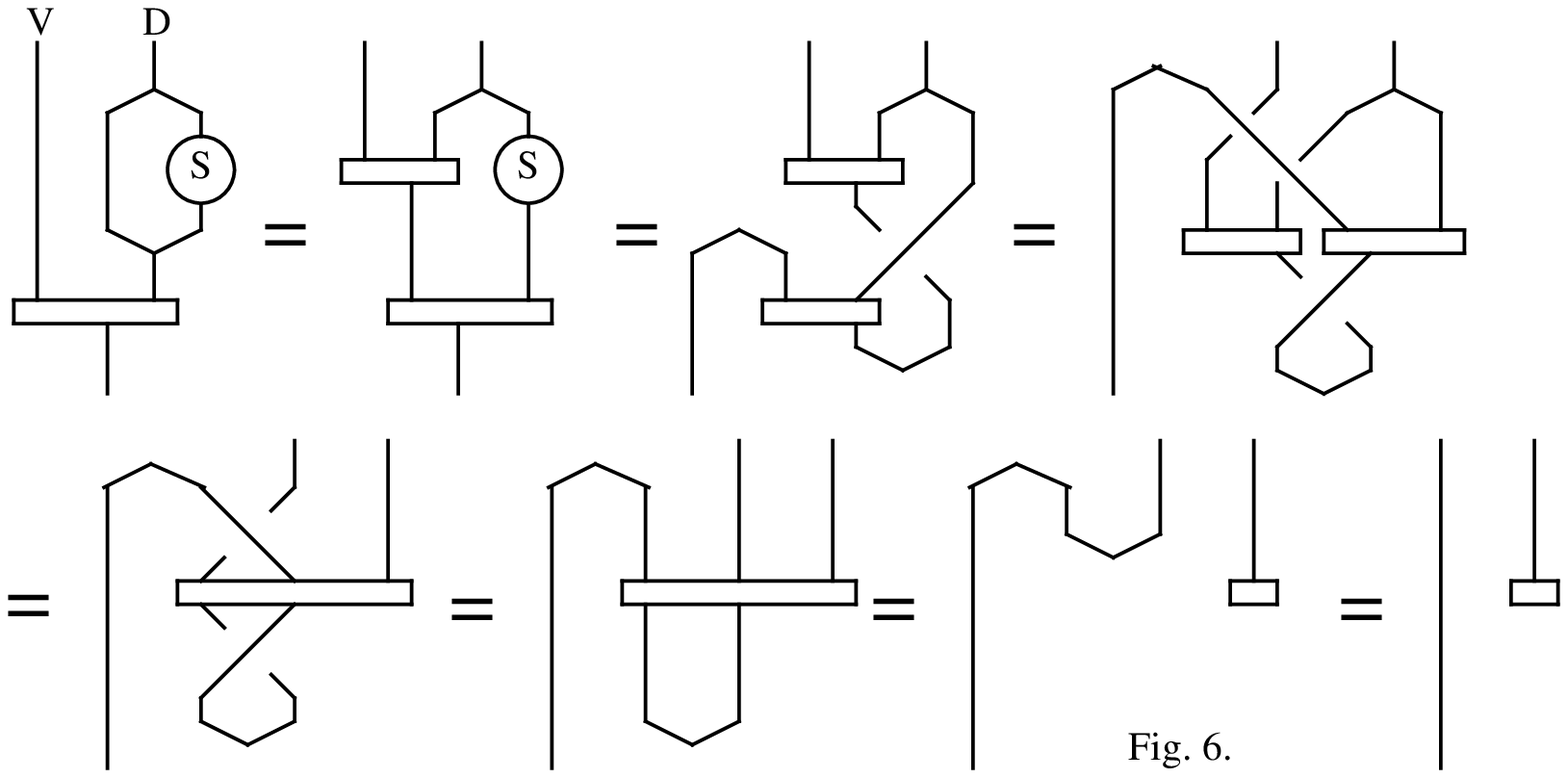}
$$
%\psboxscaled{1000}{fig6.ps}%{Fig\ 2}
Then for $\mu(S\tens I)\Delta:D\to D$ we have Fig.\ 7.
$$
\epsfig{file=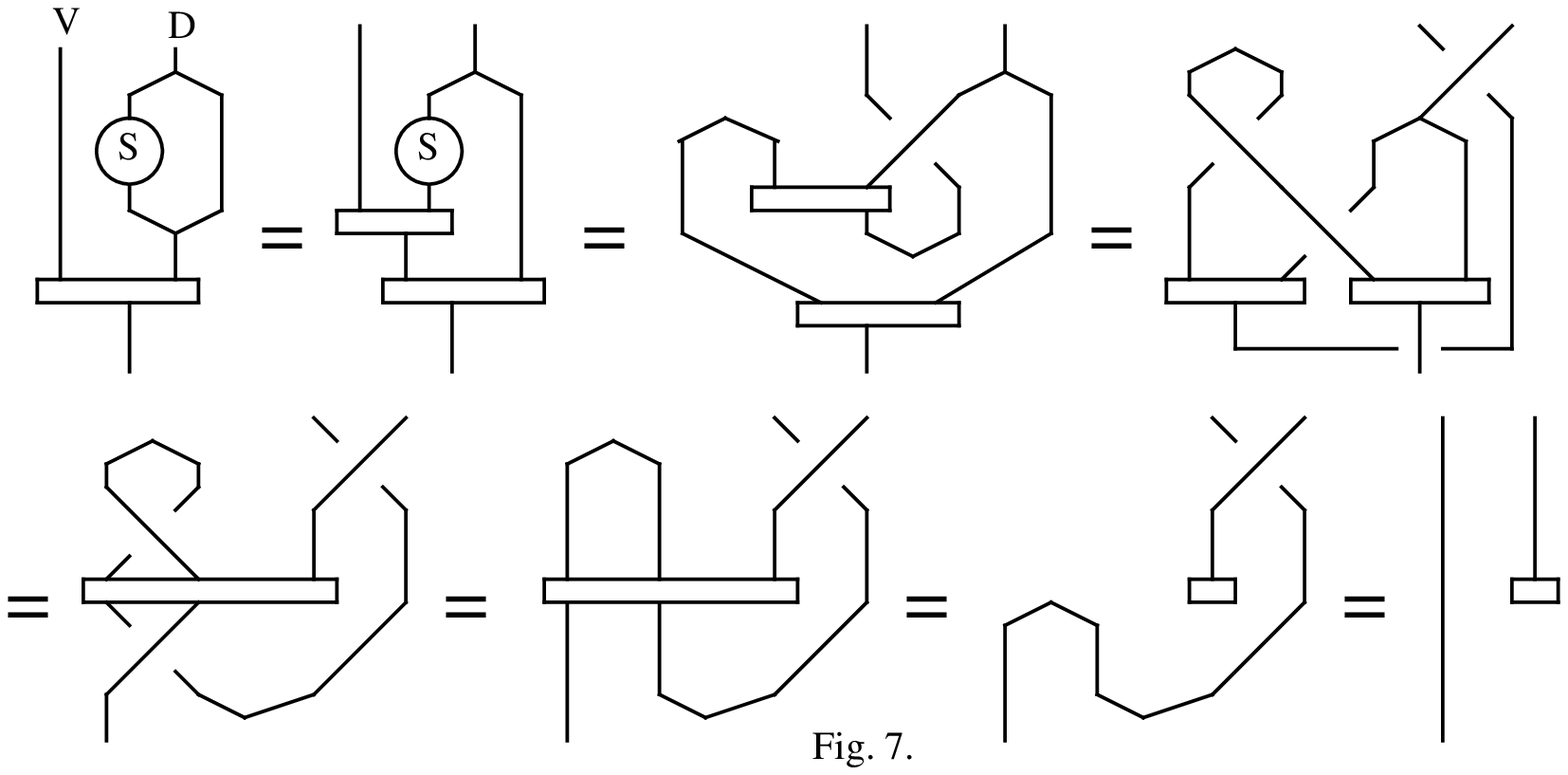}
$$
%\psboxscaled{1000}{fig7.ps}%{Fig\ 2}

\end{document}